\documentclass{article}

\usepackage{times}
\usepackage{amssymb}

\providecommand{\keywords}[1]
{
  \small	
  \textbf{\textit{Keywords---}} #1
}

\title{From the Notebooks to the Investigations\\ and Beyond\thanks{Preprint. Final version to be published in De Gruyter's \emph{SATS: Northern European Journal of Philosophy.}}}

\author{
Ruy J.G.B. de Queiroz\\
Centro de Inform\'atica\\
Universidade Federal de Pernambuco\\
Recife, Brazil\\
ruy@cin.ufpe.br}

\begin{document}

\maketitle

\begin{abstract}
The use of the open and searchable Wittgenstein's \emph{Nachlass} (The Wittgenstein Archives at the University of Bergen (WAB))  has proved instrumental in the quest for a common thread of Wittgenstein's view on the connections between meaning, use and consequences, going from the \emph{Notebooks} to later writings (including the \emph{Philosophical Investigations}) and beyond. Here we take this as the basis for a proposal for a formal counterpart of a `meaning-as-use' (dialogical/game-theoretical) semantics for the language of predicate logic. In order to further consolidate this perspective, we shall need to bring in key excerpts from Wittgenstein oeuvre (including the \emph{Nachlass}) and from those formal semanticists who advocate a different perspective on the connections between proofs and meaning. With this in mind we consider several passages from Wittgenstein's published as well as unpublished writings to build a whole picture of a formal counterpart to `meaning is use' based on the idea that explanations of consequences via `movements within language' ought to be taken as a central aspect to Wittgenstein's shift from `interpretation of symbols in a state of affairs' to `use of symbols' which underpins his `meaning is use' paradigm. As in the \emph{Investigations} ``every interpretation hangs in the air together with what it interprets, and cannot give it any support. Interpretations by themselves do not determine meaning", as well as in a remark from his transitional period (1929-30): ``Perhaps one should say that the expression ``interpretation of symbols" is misleading and one should instead say ``the use of symbols"."\footnote{Part of the material reported here was presented in a talk entitled `Explanation of Consequences via Movements within Language' at the \emph{VIII Brazilian Society for Analytic Philosophy Conference}, 22-26 July 2024 https://sites.google.com/view/sbfa-sbpha/olinda-2024\_1, 
Academia Santa Gertrudes, Olinda, Pernambuco, Brazil. Thanks to Marcos Silva for the kind invitation to give a presentation as a keynote speaker. 

We have made intensive use of Wittgenstein's \emph{Nachlass} much more than it was done in any previous publications. Thus, as a new step in a series of essays going back a few decades, it contains some overlap with the previous publications intended to bring context and self-containedness to the present manuscript, this substantially reinforces earlier arguments. (Gratitude to the editors of \emph{Wittgenstein Archives at the University of Bergen} (\emph{WAB}) who has rewarded Wittgenstein scholarship with such an immense gift as the \emph{Nachlass} in free online and searchable form!)} Significantly, we wish the present examination of the searchable \emph{Nachlass} can make a relevant step towards a formal counterpart to the ``meaning is use" dictum, while highlighting an important common thread from Wittgenstein's very early to very late writings. For this we focus here on the themes of explanation of consequences, movements within language, moving away from states of affairs, meaning vs verification.

\end{abstract}


\section{Introduction}
Supported by an extensive dive into Wittgenstein's Nachlass, 
the aim here is to consider the so-called rules of proof reduction as a formal counterpart to the explanation of the (immediate) consequences of a proposition. This contrasts markedly to a different view by the verificationist theories of meaning as put forward by Heyting, Gentzen, Dummett, Prawitz, Martin-L\"of and several others, since we suggest an approach which has more of a `pragmatist' slant to the semantics of predicate logic. Accordingly,\footnote{Here it is important to emphasise that, similarly to those who advocate the verificationist interpretation of intuitionistic logic, we are also drawing from Wittgenstein's paradigm `meaning is use', but in a way that is very different from, e.g.\ Dummett's:
\begin{quote}
 This sketch of one possible route to an account of why, within mathematics, classical logic must be abandoned in favour of intuitionistic logic obviously leans heavily upon Wittgensteinian ideas about language. (`The philosophical basis of intuitionistic logic' (1975))
 \end{quote}
 as well as from Prawitz' ``As pointed out by Dummett, this whole way of arguing with its stress on communication and the role of the language of mathematics is inspired by ideas of Wittgenstein and is very different from Brouwer's rather solipsistic view of mathematics as a languageless activity. Nevertheless, as it seems, it constitutes the best possible argument for some of Brouwer's conclusions." (`Meaning and proofs: on the conflict between classical and intuitionistic logic' (1977))  and
 \begin{quote}
 I have furthermore argued that the rejection of the platonistic theory of meaning depends, in order to be conclusive, on the development of an adequate theory of meaning along the lines suggested in the above discussion of the principles concerning meaning and use. Even if such a Wittgensteinian theory did not lead to the rejection of classical logic, it would be of great interest in itself. (Ibid.)
 \end{quote}
we take as the main ingredient of the above mentioned paradigm the \emph{application}, the \emph{purpose}, the (immediate) \emph{consequences} of the logical constant, rather than its assertability conditions. We have identified a few passages in Wittgenstein's writings where he suggests a connection between the meaning of a proposition with its method of verification, in a couple of those passages saying explicitly that the meaning of a sentence is determined by its proof, but later withdrawing his own writing and keeping `method of verification' but in a different perspective, and in this case in line with a sort of common line of thought connecting meaning and use, application, consequences.} we consider several passages from Wittgenstein's published as well as unpublished writings, in particular, The Wittgenstein Archives at the University of Bergen (WAB), to build a whole picture of a formal counterpart to `meaning is use' on the basis of the idea that explanations of consequences via `movements within language' ought to be taken as a central aspect to Wittgenstein's shift from `interpretation of symbols' to `use of symbols' which underpins his `meaning is use' paradigm. As in the \emph{Investigations}, ``every interpretation hangs in the air together with what it interprets, and cannot give it any support. Interpretations by themselves do not determine meaning", as well as in a remark from his transitional period (1929-30): 

\begin{quote}
Vielleicht mu{\ss} man sagen da{\ss} der Ausdruck ``Interpretation von Symbolen" irref\"uhrend ist und man sollte statt dessen sagen ``der Gebrauch von Symbolen". Denn ``Interpretation" klingt so als w\"urde man nun dem Wort ``rot" die Farbe rot zuordnen (wenn sie gar nicht da ist) u.s.w.. Und es entsteht wieder die Frage: Was ist der Zusammenhang zwischen Zeichen und Welt. K\"onnte ich nach etwas suchen, wenn nicht der Raum da w\"are, worin ich es suche?!

     Wo kn\"upft das Zeichen an die Welt an?
	 			 	
     Etwas suchen ist gewi{ss} ein Ausdruck der Erwartung. D.h.: Wie man sucht, dr\"uckt irgendwie aus, was man erwartet.
 			 	
     Die Idee w\"are also, da{\ss} das, was die Erwartung mit der Realit\"at gemeinsam hat, ist, da{\ss} sie sich auf einen andern Punkt im selben Raum bezieht. (Raum ganz allgemein verstanden).\footnote{which can be translated as:
     \begin{quote}
     Perhaps one should say that the expression ``interpretation of symbols" is misleading and one should instead say ``the use of symbols". For ``interpretation" sounds as if one would now assign the colour red to the word ``red" (when it is not there at all) and so on. And the question arises again: What is the connection between sign and world? Could I search for something if the space in which I am searching for it was not there?
     
     Where does the sign link up with the world?
	 			 	
     Looking for something is certainly an expression of expectation. I.e.: How one searches somehow expresses what one expects.
 	 			 	
     So the idea would be that what expectation has in common with reality is that it refers to another point in the same space. (Space understood in a general way).
     \end{quote}
Ts-208,136r, Wittgenstein Nachlass Ts-208 [based on MSS 105, 106, 107 and the first half of Ms-108] (WL)
User filtered transcription. In: Wittgenstein, Ludwig: Interactive Dynamic Presentation (IDP) of Ludwig Wittgenstein's philosophical Nachlass. Edited by the Wittgenstein Archives at the University of Bergen (WAB) under the direction of Alois Pichler. Bergen: Wittgenstein Archives at the University of Bergen 2016-. (Accessed 18 Jul 2024)
}    
\end{quote}
but also present in early writings from Ms-103 (\emph{Notebooks}), as far back as 11.9.16:
\begin{quote}
Die Art und Weise wie die Sprache bezeichnet spiegelt sich in ihrem Gebrauche wieder.\footnote{which can be translated as
\begin{quote}
The way the language describes is reflected in its use.
\end{quote}
Ms-103.51r, Wittgenstein Nachlass Ms-103 [so-called WW1 notebooks] (WL)
User filtered transcription. (WAB) (Accessed 18 Jul 2024)
}\end{quote}

In another item from the \emph{Nachlass}, this time from Ms 110 VI. \emph{Philosophische Bemerkungen} (1930--31), one finds a remark pertaining to `interpretation' of a sign:
\begin{quote}
Eine Interpretation ist immer nur eine im Gegensatz zu einer andern. Sie h\"angt sich an das Zeichen \& reiht es in ein weiteres System ein.\footnote{which can be translated as:
\begin{quote}
An interpretation is always only one in contrast to another. It attaches itself to the sign \& places it in a further system.
\end{quote}
Ms-110,288, Wittgenstein Nachlass Ms-110: VI, Philosophische Bemerkungen (WL).
User filtered transcription. In: Wittgenstein, Ludwig: Interactive Dynamic Presentation (IDP) of Ludwig Wittgenstein's philosophical Nachlass. Edited by the Wittgenstein Archives at the University of Bergen (WAB) under the direction of Alois Pichler. Bergen: Wittgenstein Archives at the University of Bergen 2016-. (Accessed 29 Dec 2023)
}
\end{quote}

Since this might sound like it is leaning towards a `pragmatist' perspective on language and meaning, it is relevant to quote from Wittgenstein himself on the question of whether he considered his account of the connections between meaning, use and purpose as a pragmatist one:
\begin{quote}
Hier sieht man den Zugang zu der pragmatistischen Auffassung von Wahr \& Falsch. Der Satz ist solange wahr solang er sich als n\"utzlich erweist.

Jeder Satz den wir im gew\"ohnlichen Leben \"au{\ss}ern scheint den Charakter einer Hypothese zu haben.
	 			 	
Die Hypothese ist ein logisches Gebilde. D.h.\ ein besonderes Symbol wof\"ur bestimmte || gewisse Regeln der Darstellung gelten.
	 			 	
Das Reden von Sinnesdaten \& der unmittelbaren Erfahrung, hat den Sinn, da\ss\ wir eine nicht-hypothetische Darstellung suchen.
 			 	
Nun scheint es aber da\ss\ die Darstellung \"uberhaupt ihren Wert verliert wenn man das hypothetische Element in ihr fallen l\"a{\ss}t, weil dann der Satz nicht mehr auf die Zukunft deutet sondern quasi selbstzufrieden ist \& daher wertlos.
 			 	
Die Erfahrung sagt gleichsam ,,sch\"on ist es auch anderswo \& hier bin ich sowieso". Und mit dem Perspektiv der Erwartung schauen wir in die Zukunft.
 	 			 	
Es hat keinen Sinn von S\"atzen zu reden die als Instrumente keinen Wert haben.
 	 			 	
Der Sinn eines Satzes ist sein Zweck.\footnote{which can be translated as:
\begin{quote}
Here you can see the approach to the pragmatist concept of true and false. The sentence is true as long as it proves to be useful.

Every proposition we utter in ordinary life seems to have the character of a hypothesis.
 	 			 	
The hypothesis is a logical construct. I.e.\ a particular symbol to which certain rules of representation apply.
	 			 	
Speaking of sense data \& immediate experience, has the sense that we seek a non-hypothetical representation.
 	 			 	
But now it seems that the representation loses its value at all if one drops the hypothetical element in it, because then the proposition no longer points to the future but is quasi self-satisfied \& therefore worthless.
 	 			 	
The experience says, as it were, `it's also nice elsewhere \& I'm here anyway'. And with the perspective of expectation, we look to the future.
		 	
There is no point in talking about sentences that have no value as instruments.
		 	
The meaning of a sentence is its purpose.
\end{quote}
Ms-107,248-249 (Ms-107: III, \emph{Philosophische Betrachtungen}, Sept 1929 to 5 Jan 1930) (Access 18 July 2024)
}
\end{quote}
In order to answer a possible question as to whether he is a pragmatist, he comes up with:
\begin{quote}
Wie aber, wenn die Religion lehrt, die Seele k\"onne bestehen, wenn der Leib zerfallen ist? Verstehe ich, was sie lehrt? Freilich verstehe ich's: ? ? Ich || ich kann mir dabei manches vorstellen.

(Man hat ja auch Bilder von diesen Dingen gemalt. Und warum sollte so ein Bild nur die unvollkommene Wiedergabe des ausgesprochenen Gedankens sein? Warum soll es nicht den gleichen Dienst tun, wie der Satz? || , wie das, was wir sagen?)

     Und auf den Dienst kommt es an.
 	 			 	
     Aber bist Du kein Pragmatiker? Nein. Denn ich sage nicht, der Satz sei wahr, der n\"utzlich ist.
     
     Der Nutzen, d.h., Gebrauch, gibt dem Satz seinen besondern Sinn, das Sprachspiel gibt ihm ihn.
     
     Und insofern als eine Regel oft so gegeben wird, da\ss\ sie sich n\"utzlich erweist, \& mathematische S\"atze mit Regeln || ihrer Natur || ihrem Wesen nach mit Regeln verwandt sind, spiegelt sich in mathematischen
Wahrheiten N\"utzlichkeit.\footnote{which can be translated as:
\begin{quote}
But what if religion teaches that the soul can exist when the body has decayed? Do I understand what she teaches? Of course I understand it:  I || I can imagine a lot of things.

(People have also painted pictures of these things. And why should such a picture only be an imperfect representation of the thought expressed? Why shouldn't it do the same job as the sentence? || , as what we say?)
     
     And the service is what matters.

      But aren't you a pragmatist? No. For I do not say that the sentence that is useful is true.
      
      Utility, i.e. use, gives the sentence its special meaning; the language game gives it this.
      
      And insofar as a rule is often given in such a way that it proves useful, \& mathematical sentences with rules || their nature || are essentially related to rules, usefulness is reflected in mathematical truths.
\end{quote}
Ms-131,69-71 (1946),  Ts-245,185 (1945--47), and Ts 229,253 (1947) Wittgenstein Nachlass Ms-131 (WL). (Accessed 01 October 2024)
}
\end{quote}

As it happens, one can find traces of such a connection between the meaning of a sign and its use, application, employment, very early in his writings, such as in:
\begin{quote}
Namen kennzeichnen die Gemeinsamkeit einer Form und eines Inhalts. -- Sie kennzeichnen erst mit ihrer syntaktischen Verwendung zusammen eine bestimmte logische Form.\footnote{which can be translated as:
\begin{quote}
Names characterise the commonality of a form and a content. - They only characterise a certain logical form together with their syntactic use.
\end{quote}
Ms-102,116r, 30.5.15 (Also in Ms-104, 3'253, \emph{Prototractatus}) (Accessed 01 October 2024)
}
\end{quote}
A similar remark appears in the \emph{Prototractatus}:
\begin{quote}	 			 	
3'253 Zeichen kennzeichnen die Gemeinsamkeit einer Form und eines Inhalts. -- Sie bestimmen erst mit ihrer syntaktischen Verwendung zusammen eine logische Form.\footnote{which can be translated as:
\begin{quote}
3'253  Signs characterise the commonality of a form and a content. - They only determine a logical form together with their syntactic use.
\end{quote}
Ms-104,58 [Ms-104: Logisch-Philosophische Abhandlung [so-called Prototractatus] (Bodl, MS. German d. 7)] (Accessed 01 October 2024)
}
\end{quote}
Roughly two weeks later, Wittgenstein writes:
\begin{quote}
Wir sind uns also dar\"uber klar geworden da\ss\ Namen f\"ur die verschiedensten Formen stehen, und stehen d\"urfen, und da\ss\ nun erst die syntaktische Verwendung / Anwendung die darzustellende Form charakterisiert. Was ist nun die syntaktische Anwendung von Namen einfacher Gegenst\"ande?\footnote{which can be translated as:
\begin{quote}
We have thus become aware that names stand, and may stand, for the most diverse forms, and that only the syntactic use characterises the form to be represented. What is the syntactic use of names of simple objects?
\end{quote}
Ms-102,141r, 14.6.15 (Accessed 01 October 2024)
}
\end{quote}
Just two days after this, here comes another remark on the `syntaktische Anwendung' of a name:
\begin{quote}
Das hei{\ss}t die syntaktische Verwendung von || der Namen charakterisiert vollst\"andig die Form der zusammengesetzten Gegenst\"ande welche sie bezeichnen.\footnote{which can be translated as:
\begin{quote}
This means that the syntactic use of || of the names completely characterises the form of the compound objects they denote.
\end{quote}
Ms-102,149r, 16.6.15
}
\end{quote}
And then, again two days later it comes a remark which does not apply just to names but to propositional signs:
\begin{quote}
Sagen wir dann von einem Punkt in diesem || jenem Fleck da\ss\ er rechts von der Linie liege dann folgt dieser Satz aus dem fr\"uheren und wenn unendlich viele Punkte in dem Flecken liegen dann folgen unendlich viele S\"atze verschiedenen Inhalts logisch aus jenem ersten! Und dies zeigt schon da\ss\ er tats\"achlich selbst unendlich komplex war. N\"amlich nicht das Satzzeichen allein wohl aber mit seiner syntaktischen Verwendung.\footnote{which can be translated as:
\begin{quote}
If we then say of a point in this || that spot that it lies to the right of the line then this proposition follows from the earlier one and if infinitely many points lie in the spot then infinitely many propositions of different content follow logically from that first one! And this already shows that he himself was actually infinitely complex. Not the propositional sign alone, but its syntactical use.
\end{quote}
Ms-102,162r, 18.6.15
}
\end{quote}
this is followed almost immediately by a similar remarks concerning not only names and propositional signs, but to the syntactic use of each component of a sentence in order for it to have a meaning:
\begin{quote}
Man k\"onnte die Bestimmtheit auch so fordern!: Wenn ein Satz Sinn haben soll so mu\ss\ vorerst die syntaktische Verwendung jedes seiner Teile festgelegt sein. -- Man kann z.B.\ nicht erst nachtr\"aglich draufkommen da\ss\ ein Satz aus ihm folgt. Sondern z.B.\ welche S\"atze aus einem Satz folgen mu\ss\ vollkommen feststehen ehe dieser Satz einen Sinn haben kann!\footnote{which can be translated as:
\begin{quote}
One could also demand certainty like this: If a sentence is to have meaning, the syntactical use of each of its parts must first be determined. -- For example, one cannot discover afterwards that a sentence follows from it. For example, which sentences follow from a sentence must be completely determined before this sentence can have a meaning!
\end{quote}
Ms-102,164r, 18.6.15 (Also Ms-104,60, 3'20103 \emph{Prototractatus})
}
\end{quote}
In the manuscript Ms-104, the so-called \emph{Prototractatus} one finds:
\begin{quote}
3'11 Das Satzzeichen ist eine Projektion seines || eines || des Gedankens.
 	 			 	
3'12 Die Projektionsmethode ist die Art und Weise der Anwendung des Satzzeichens.
	 			 	
3'13 Die Anwendung des Satzzeichens ist das Denken seines Sinns.\footnote{which can be translated as:
\begin{quote}
3'11 The propositional sign is a projection of its || of a || of thought.
 	 			 	
3'12 The method of projection is the way in which the propositional sign is used.
			 	
3'13 The application of the propositional sign is the thinking of its meaning.
\end{quote}
Ms-104,5, 
}
\end{quote}

Going forward, we come to the \emph{Tractatus} where Wittgenstein makes use of the words Anwendung/Verwendung to associate meaning and application/employment/use:
\begin{quote}
3.262 Was in den Zeichen nicht zum Ausdruck kommt, das zeigt ihre Anwendung. Was die Zeichen verschlucken, das spricht ihre Anwendung aus.\footnote{which was translated by Ogden as:
\begin{quote}
3.262 What does not get expressed in the sign is shown by its application. What the signs conceal, their application declares.
\end{quote}
and translated by Pears \& McGuinness as:
\begin{quote}
3.262 What signs fail to express, their application shows. What signs slur over, their application says clearly.
\end{quote}
but can also be translated as:
\begin{quote}
3.262 What is not expressed in the signs is shown by their use. What the signs swallow up is expressed by their use.
\end{quote}
}
\end{quote}
Further down comes another instance of this:
\begin{quote}
3.326 Um das Symbol am Zeichen zu erkennen, muss man auf den sinnvollen Gebrauch achten.

3.327 Das Zeichen bestimmt erst mit seiner logisch-syntaktischen Verwendung zusammen eine logische Form.\footnote{which was translated by Ogden as:
\begin{quote}
3.326 In order to recognize the symbol in the sign we must consider the significant use.

3.327 The sign determines a logical form only together with its logical syntactic application.
\end{quote}
and translated by Pears \& McGuinness as:
\begin{quote}
3.326 In order to recognize a symbol by its sign we must observe how it is used with a sense.

3.327 A sign does not determine a logical form unless it is taken together with its logico-syntactical employment.
\end{quote}
but can also be translated as:
\begin{quote}
3.327 The sign only determines a logical form together with its logical-syntactical use.
\end{quote}
}
\end{quote}

Furthermore, one can find remarks in the \emph{Tractatus} pertaining to the connection between meaning and purpose, starting with this:
\begin{quote}
3.341 Das Wesentliche am Satz ist also das, was allen S\"atzen, welche den gleichen Sinn ausdr\"ucken k\"onnen, gemeinsam ist.	

Und ebenso ist allgemein das Wesentliche am Symbol das, was alle Symbole, die denselben Zweck erf\"ullen k\"onnen, gemeinsam haben..\footnote{which has been translated as:
\begin{quote}
3.341 The essential in a proposition is therefore that which is common to all propositions which can express the same sense.	

And in the same way in general the essential in a symbol is that which all symbols which can fulfill the same purpose have in common. (Ogden)
\end{quote}
\begin{quote}
3.341 So what is essential in a proposition is what all propositions that can express the same sense have in common.

And similarly, in general, what is essential in a symbol is what all symbols that can serve the same purpose have in common. (Pears \& McGuinness)
\end{quote}
which can also be translated as:
\begin{quote}
3.341 The essence of the sentence is therefore what all sentences that can express the same meaning have in common.

And likewise, the essence of the symbol in general is what all symbols that can fulfil the same purpose have in common.
\end{quote}
}
\end{quote}
As a matter of fact, this remark is already in the \emph{Notebooks} in a slightly different way:
\begin{quote} 			 	
Es ist klar da\ss\ Zeichen, die denselben Zweck erf\"ullen logisch identisch sind. Das rein Logische ist eben das was alle diese leisten k\"onnen.\footnote{which can be translated as:
\begin{quote}
It is clear that signs that fulfil the same purpose are logically identical. The purely logical is precisely what all of them can do.
\end{quote}
Ms-102,75r, 23.4.15. (Ms-102 [so-called WW1 notebooks] (WL)) (Access 23/09/2024)
}
\end{quote}
Concerning the rules of logical syntax and what the sign signifies, there is a remark very much near the abovementioned one:
\begin{quote}
3.344 Das, was am Symbol bezeichnet, ist das Gemeinsame aller jener Symbole, durch die das erste den Regeln der logischen Syntax zufolge ersetzt werden kann.\footnote{which has been translated as:
\begin{quote}
3.344 What signifies in the symbol is what is common to all those symbols by which it can be replaced according to the rules of logical syntax. (Ogden)
\end{quote}
\begin{quote}
3.344 What signifies in a symbol is what is common to all the symbols that the rules of logical syntax allow us to substitute for it. (Pears \& McGuinness)
\end{quote}
which can also be translated as:
\begin{quote}
3.344 That which is denoted by the symbol is common to all those symbols by which the first can be replaced according to the rules of logical syntax..
\end{quote}
}
\end{quote}
Further down, one finds a direct link between sameness of purpose and sameness of meaning:
\begin{quote}
5.47321 Zeichen, die Einen Zweck erf\"ullen, sind logisch \"aquivalent, Zeichen, die keinen Zweck erf\"ullen, logisch bedeutungslos.\footnote{which has been translated as:
\begin{quote}
5.47321 Signs which serve \emph{one} purpose are logically equivalent, signs which serve \emph{no} purpose are logically meaningless. (Ogden)
\end{quote}
\begin{quote}
5.47321 Signs that serve \emph{one} purpose are logically equivalent, and signs that serve \emph{none} are logically meaningless. (Pears \& McGuinness)
\end{quote}
which can also be translated as:
\begin{quote}
5.47321 Signs that fulfil a purpose are logically equivalent, signs that do not fulfil a purpose are logically meaningless.
\end{quote}
}
\end{quote}
All the way until the very last writings one can find the suggestion that the meaning of a word, a term or a sentence is to be found in the use, application, purpose, etc., be it in the form `Gebrauch', or in German words which are more frequent in the very early writings, including the \emph{Notebooks 1914--16}. Here are some passages from manuscripts written in 1937--38 and 1945:
\begin{quote}
Die Anwendung des Satzes ist nicht die, die ein solches Vorstellen fordert. Immer wieder m\"ochte man sich den Sinn eines Satzes, also seine Verwendung (seinen Nutzen) in einem seelischen Zustand des Redenden oder H\"orenden konzentriert denken. Man denkt nicht, da\ss\ man mit den Worten rechnet, operiert, f\"ur sie mit der Zeit dies oder jenes Bild substituiert. || sie mit der Zeit in dies oder jenes Bild \"uberf\"uhrt. Sondern der -- ihr Sinn, d.i.\ aber ihr Zweck, soll in einer Art Bild liegen, das sie im Geist des Sprechers erzeugen. Es ist so -- ganz so als glaubte -- glaube man, da\ss\ etwa einer schriftlichen Anweisung, auf eine Kuh, das gefolgt werden soll, --| die mir Einer -- irgend jemand ausfolgen soll, immer die Vorstellung von einer Kuh folgen m\"usse || die schriftliche Anweisung, auf eine Kuh, die mir Einer -- irgend jemand ausfolgen soll, immer von einer Vorstellung einer Kuh begleitet sein -- werden m\"usse, wenn diese Anweisung nicht ihren Sinn verlieren soll. || damit diese Anweisung nicht ihren Sinn verliere.\footnote{which can be translated as:
\begin{quote}
The application of the sentence is not the one that demands such an imagination. Again and again one would like to think of the meaning of a sentence, i.e.\ its use (its utility) in a mental state of the speaker or hearer in a concentrated way. One does not think that one is calculating, operating with the words, substituting this or that image for them over time. -- One does not think that one transforms them into this or that image over time. Rather, their meaning, i.e.\ their purpose, should lie in a kind of image that they create in the mind of the speaker. It is as if one believed that, for example, a written instruction to follow a cow, which someone is to follow, must always be followed by an image of a cow, if this instruction is not to lose its meaning. -- so that this instruction does not lose its meaning.
\end{quote}
Ms-116,83 [Ms-116: XII, Philosophische Bemerkungen (WL)] (1937-38, 1945) (Accessed 06 October 2024)
}
\end{quote}
The pattern is carried along all the way into very late writings such as those from
1950--1951:
\begin{quote}
Als ich sagte es sei eine Unbestimmtheit in der Anwendung, meinte ich nicht, ich wisse nicht recht, wann ich die \"Au{\ss}erung machen solle (wie es etwa w\"are, wenn ich nicht recht || gut Deutsch verst\"unde).
 			 	
     Man mu\ss\ -- darf eben nicht vergessen, welche Verbindungen gemacht werden, wenn wir lernen Ausdr\"ucke wie ``Ich \"argere mich" zu gebrauchen.
	 			 	
     Und denke nicht an ein Erraten der richtigen Bedeutung durch das Kind, denn, ob es sie richtig erraten hat, mu\ss\ sich doch wieder in seiner Verwendung der Worte zeigen.\footnote{which can be translated as:
\begin{quote}
When I said it was an indeterminacy in the application,
I didn't mean that I didn't really know when to make the utterance (as it would be if I didn't understand German very well).
	 			 	
     We must not forget what connections are made when we learn to use expressions like `Ich \"argere mich'.
	 			 	
     And don't think about the child guessing the correct meaning, because whether he has guessed it correctly must be shown in his use of the words.
\end{quote}
Ms-173,47r [Ms-173 [so-called Notebook no.5; source for Ms-176] (WL)] (between 24 March and 12 April 1950) (Accessed 04 October 2024)
}

Und wie kann es unsinnig sein, zu sagen ``Es gibt Menschen, welche sehen", wenn es nicht unsinnig ist, zu sagen, es gibt Menschen, welche blind sind?

     Aber der Sinn des Satzes ``Es gibt Menschen, welche sehen" d.h.\ seine m\"ogliche Verwendung ist jedenfalls nicht sogleich klar.\footnote{which can be translated as:
\begin{quote}
And how can it be nonsensical to say `There are people who see' if it is not nonsensical to say that there are people who are blind?

     But the meaning of the sentence ?There are people who see?, i.e.\ its possible use, is not immediately clear.
\end{quote}
Ms-173,97r [Ms-173 [so-called Notebook no.5; source for Ms-176] (WL)] (between 24 March and 12 April 1950) (Accessed 04 October 2024)
}
     
\end{quote}
And again:
\begin{quote}
Wir sagen: wenn -- damit -- wenn das Kind die Sprache, --  \& also ihre Anwendung -- beherrscht, mu\ss\ es die Bedeutungen der Worte wissen. Es mu\ss\ z.B.\ einem wei{\ss}en, schwarzen, roten, oder
grauen -- blauen Dinge seinen Farbnamen, in der Abwesenheit jedes Zweifels, beilegen k\"onnen.

Ja, hier vermi{\ss}t auch niemand den Zweifel; wundert sich niemand, da\ss\ wir die Bedeutung der Worte nicht nur vermuten.\footnote{which can be translated as:
\begin{quote}
We say: if -- so -- if the child masters the language, -- - \& thus its use - it must know the meanings of the words. He must, for example, be able to give a white, black, red or grey thing its colour name in the absence of any doubt.
	 			 	
Yes, here no one misses the doubt; no one is surprised that we do not merely surmise the meaning of the words.
\end{quote}
Ms-176,46v [Ms-176 [so-called Notebook no.6.1] (WL)] (1951) (Accessed 04 October 2024)
}

Das Kind, welches die Verwendung dieses -- des Wortes beherrscht, ist nicht `sicher, diese Farbe hei{\ss}e in seiner Sprache so'. Man kann auch nicht von ihm sagen, es lerne, wenn es sprechen lernt, da\ss\ die Farbe auf Deutsch `rot' || so hei{\ss}t, indem es sprechen lernt; \& -- oder auch, || : es wisse dies, wenn es den Gebrauch des Worts erlernt hat.\footnote{which can be translated as:
\begin{quote}
The child who has mastered the use of this word is not `sure that this colour is called that in his language'. Nor can he be said to learn, when he learns to speak, that the colour is called `red' in German by learning to speak; \& -- or even, -- : he knows this when he has learned the use of the word.
\end{quote}
Ms-176,53v [Ms-176 [so-called Notebook no.6.1] (WL)] (1951) (Accessed 04 October 2024)
}
\end{quote}

\section{Moving away from `states of affairs'}

Already in the manuscripts from 1929, right after his return to writing philosophy, there is a clear attempt to reassess the fundamental role of the notion of `state of affairs' which served as the basis for what is considered his `picture theory'. In an explicit recognition of the restricted view of Augustine about language and reality, especially with respect to names and their meanings, a quote from an English translation of a prewar version of the \emph{Investigations} of 1938-39 says:
\begin{quote}
It is interesting to compare the variety of the instruments of our language and of their applications -- the ways they are applied -- their various uses --  the variety of the parts of speech and of the kinds of -- kinds of words \& of sentences -- with what logicians have said about the structure of our language. (And the author of the Tractatus Logico-philosophicus as well -- including the author of Tract. Log.phil.)

This is connected with the view that the -- fact that we think that the learning of the language consists in naming objects; namely -- viz.\ human beings, forms -- shapes, colours, pains -- aches, moods, numbers, etc..- As we have said, - naming is something like affixing a nameplate to -- putting -- fastening a label to a thing. One may call this -- And this one might call the -- a preparation for the use of a word. But for what is it a preparation?\footnote{Ts-226,16, Ts-226 [English translation of part of the Investigations prewar version] (WL) (1938-39) (Accessed 30 September 2024)
}
\end{quote}

As a matter of fact, a significant amount of the introductory material of those first writings of this period brings in a few elements of Euclidean space in order to unveil 
the limitations of the idea of a coordinate system to make the bridge between language and reality. In his own words:
\begin{quote}
Als ich die Sprache ersann || konstruierte die sich bei der Darstellung des Sachverhaltes im Raum eines Koordinatensystems bedient da habe ich doch damit einen Bestandteil in die Sprache eingef\"uhrt dessen sie sich sonst nicht bedient. Dieses Mittel ist gewi\ss\ erlaubt. Und es zeigt den Zusammenhang zwischen Sprache \& Realit\"at. Das geschriebene Zeichen ohne das Koordinatensystem ist sinnlos. Mu\ss\ nun nicht etwas \"ahnliches im Falle || zur Darstellung der Farben verwendet werden?\footnote{which can be translated as:
\begin{quote}
When I constructed the language that uses a coordinate system to represent the state of affairs in space, I introduced a component into the language that it does not otherwise use. This means is certainly permitted. And it shows the connection between language and reality. The written symbol without the coordinate system is meaningless. Shouldn't something similar be used in the case of || to represent the colours?
\end{quote}
Ms-107,280 [Ms-107: III, Philosophische Betrachtungen] (1929-30) (Accessed 29 September 2024)
}
\end{quote}
Further to this acknowledgment of his own misconception, there is a reference to `those who only attribute reality to things':
\begin{quote}
Da\ss\ uns nichts auff\"allt, wenn wir uns umsehen, im Raum herumgehen, unseren eigenen K\"orper f\"uhlen etc. etc. das zeigt, wie nat\"urlich uns eben diese Dinge sind. Wir nehmen nicht wahr, da\ss\ wir den Raum perspektivisch sehen oder da\ss\ das Gesichtsbild gegen den Rand zu in irgend einem Sinne verschwommen ist. Es f\"allt uns nie auf und kann uns nie auffallen, weil es die Art der Wahrnehmung ist. Wir denken nie dar\"uber nach, und es ist unm\"oglich, weil es zu der Form unserer Welt keinen Gegensatz gibt.

 Ich wollte sagen, es ist merkw\"urdig, da\ss\ die, die nur den Dingen, nicht unseren Vorstellungen, Realit\"at zuschreiben, sich in der Vorstellungswelt so selbstverst\"andlich bewegen und sich nie aus ihr heraussehnen.
 
     D.h., wie selbstverst\"andlich ist doch das Gegebene. Es m\"u{\ss}te mit allen Teufeln zugehen, wenn das das kleine, aus einem schiefen Winkel aufgenommene Bildchen w\"are.
     
     Dieses Selbstverst\"andliche, das Leben, soll etwas Zuf\"alliges, Nebens\"achliches sein; dagegen etwas wor\"uber ich mir normalerweise nie den Kopf zerbreche, das Eigentliche!\footnote{which can be translated as:
\begin{quote}
The fact that we don't notice anything when we look around, walk around the room, feel our own body, etc., etc., shows how natural these things are to us. We do not realise that we see the room in perspective or that the facial image is blurred towards the edge in any sense. We never notice it and can never notice it because it is the way we perceive. We never think about it, and it is impossible because there is no contrast to the form of our world.

 I wanted to say that it is strange that those who ascribe reality only to things, not to our ideas, move so naturally in the world of ideas and never long to leave it.
 
     In other words, how self-evident the given is. It would have to be a devil's bargain if it were the little picture taken from an oblique angle.
     
     This self-evident thing, life, is supposed to be something accidental, incidental; on the other hand, something that I normally never worry about, the real thing!
\end{quote}
Ts-209,17 [Ts-209: Philosophische Bemerkungen [Cuttings from Ts-208]] (1929) (Accessed 29 September 2024)
}
\end{quote}
As a followup, one sees the sign of a move towards a series of references to a view which purports to replace `state of affairs' as an unsatisfactory account of the bridge between language and reality by the introduction of a rather peculiar connection between meaning and verification:
\begin{quote}
Immer wieder ist es der Versuch, die Welt in der Sprache abzugrenzen und hervorzuheben -- was aber nicht geht. Die Selbstverst\"andlichkeit der Welt dr\"uckt sich eben darin aus, da\ss\ die Sprache nur sie bedeutet, und nur sie bedeuten kann.

     Denn, da die Sprache die Art ihres Bedeutens erst von ihrer Bedeutung, von der Welt, erh\"alt, so ist keine Sprache denkbar, die nicht diese Welt darstellt.
 	 			 	
     Wenn die Welt der Daten zeitlos ist, wie kann man dann \"uberhaupt \"uber sie reden?
 	 			 	
     Der Strom des Lebens, oder der Strom der Welt flie{\ss}t dahin, und unsere S\"atze werden, sozusagen, nur in Augenblicken verifiziert.
     
     Unsere S\"atze werden nur von der Gegenwart verifiziert.\footnote{which can be translated as:\begin{quote}
Again and again it is the attempt to delimit and emphasise the world in language - but this is not possible. The self-evident nature of the world is expressed in the fact that language only means it and can only mean it.

     For, since language receives the nature of its meaning only from its meaning, from the world, no language is conceivable that does not represent this world.
 	 			 	
     If the world of data is timeless, how can we talk about it at all?
	 			 	
     The stream of life, or the stream of the world, flows along, and our sentences are, so to speak, only verified in moments.
     
     Our sentences are only verified by the present.
\end{quote}
Ibid.
}
\end{quote}
Several moments in that 1929-1930 period witness the insistence on `method of verification' as an essential component of meaning:
\begin{quote}
Die Verifikation ist nicht ein Anzeichen der Wahrheit sondern der Sinn des Satzes.\footnote{which can be translated as:
\begin{quote}
Verification is not an indication of the truth but the meaning of the sentence.
\end{quote}
Ms 107,143
}
Wenn ich sage ,,in dieser Klasse fehlen jeden Tag durchschnittlich 5 Sch\"uler", was hei{\ss}t das. Wie wird es verifiziert; denn ,,wie ein Satz verifiziert wird, das sagt er".\footnote{which can be translated as:
\begin{quote}
If I say ``in this class an average of 5 students are absent every day", what does that mean? How is it verified; because ``how a sentence is verified is what it says". 
\end{quote}
Ms 107,175
}
Wie k\"onnten wir aber einen solchen Satz wissen? ? Wie wird er verifiziert?! Was dem, was wir meinen, wirklich entspricht ist (glaube ich) gar kein Satz sondern der Schlu\ss\ von $\phi$x auf $\psi$x, wenn dieser Schlu\ss\ gestattet ist ? aber der wird nicht durch einen Satz ausgedr\"uckt.\footnote{which can be translated as:
\begin{quote}
But how could we know such a sentence? - How is it verified? What really corresponds to what we mean is (I think) not a proposition at all but the conclusion from $\phi$x to $\psi$x, if this conclusion is allowed - but this is not expressed by a proposition.
\end{quote}
Ms 107,184
}
Wir k\"onnen unser altes Prinzip auf die S\"atze, die eine Wahrscheinlichkeit aussagen, anwenden \& sagen da\ss\ wir ihren Sinn erkennen werden wenn wir wissen || bedenken wie sie verifiziert werden || was sie verifiziert.\footnote{which can be translated as:
\begin{quote}
We can apply our old principle to the propositions that state a probability \& say that we will recognise their meaning when we know || consider how they are verified || what verifies them.
\end{quote}
Ms 108,49 [Ms-108: IV, Philosophische Bemerkungen (WL)] (1930) (Accessed 29 September 2024)
}
Ich will wissen was es hei{\ss}t, einen Satz zu verstehen: Wie verifiziert man denn diese Aussage?\footnote{which can be translated as:
\begin{quote}
I want to know what it means to understand a sentence: How do you verify this statement?
\end{quote}
Ms 108,63 [Ms-108: IV, Philosophische Bemerkungen (WL)] (1930) (Accessed 29 September 2024)
}
Die Verifikation ist nicht ein blo{\ss}es Anzeichen der Wahrheit, sondern sie bestimmt den Sinn des Satzes.\footnote{which can be translated as:
\begin{quote}
Verification is not a mere indication of truth, but rather it determines the meaning of the sentence.
\end{quote}
Ms 113,142v [Ms-113: IX, Philosophische Grammatik] (1930) (Accessed 30 September 2024)
}
Welche S\"atze aus ihm folgen \& aus welchen S\"atzen er folgt, das macht seinen Sinn aus. Daher auch die
Frage nach seiner Verifikation eine Frage nach seinem Sinn ist.\footnote{which can be translated as:
\begin{quote}
Which sentences follow from it \& from which sentences it follows, that makes its sense. Hence also the
question about its verification is a question about its meaning.
\end{quote}
Ms 113,42r [Ms-113: IX, Philosophische Grammatik] (1930) (Accessed 30 September 2024)
}
\end{quote}
In the writings about the foundations of mathematics, during the years 1931-1932, there appears a remark on `proofs vs meaning' of mathematical propositions which will be later withdrawn:
\begin{quote}
Die Arithmetik ist, wie gesagt, auch ohne eine Regel. Wie w\"are es, wenn man au{\ss}er den Multiplikationsregeln noch ,,25 $\times$ 25 = 625" als Regel festsetzen wollte? (Ich sage nicht ,,25 $\times$ 25 = 624"!) ? 25 $\times$ 25 = 625 hat nur Sinn, wenn die Art der Rechnung || Ausrechnung bekannt ist, die zu dieser Gleichung geh\"ort, \& hat nur Sinn in Bezug auf diese Rechnung. A hat nur Sinn mit Bezug auf die Art der Ausrechnung von A. Denn die erste Frage w\"are hier eben: ist das eine Bestimmung || Festsetzung, oder ein errechneter Satz? Denn ist 25 $\times$ 25 = 625 eine Festsetzung (Grundregel), dann bedeutet das Multiplikationszeichen jetzt etwas anders, als es z.B.\ in Wirklichkeit bedeutet. (D.h.\ wir haben es mit einer anderen Rechnungsart zu tun.) Und ist A eine Festsetzung dann definiert das die Addition anders, als wenn es ein errechneter Satz ist. Denn die Festsetzung ist ja dann eine Erkl\"arung des Additionszeichens \& die Rechenregeln, die A auszurechnen erlauben, eine andere Erkl\"arung desselben Zeichens. Ich darf hier nicht vergessen da\ss\ $\alpha,\beta,\gamma$ nicht der Beweis von A ist sondern nur die Form des Beweises, oder des Bewiesenen, ist; $\alpha,\beta,\gamma$ definiert also A.

Darum kann ich nur sagen ,,25 $\times$ 25 = 625 wird bewiesen", wenn die Beweismethode fixiert ist, unabh\"angig von dem speziellen Beweis. Denn diese Methode bestimmt erst die Bedeutung von ,,$\xi\times\eta$", also, was bewiesen
wird. Insofern geh\"ort also die Form aa : b = c zur Beweismethode, die den Sinn von c. erkl\"art. Etwas anderes ist dann die Frage, ob ich richtig gerechnet habe. -- Und so geh\"ort $\alpha,\beta,\gamma$ zur Beweismethode die den Sinn des Satzes A erkl\"art.\footnote{which can be translated as:
\begin{quote}
As I said, arithmetic is also without a rule. What would it be like if, in addition to the multiplication rules, you wanted to define ``25 $\times$ 25 = 625" as a rule? (I am not saying ``25 $\times$ 25 = 624"!) - 25 $\times$ 25 = 625 only makes sense if the type of arithmetic || calculation that belongs to this equation is known, \& only makes sense in relation to this calculation. A only makes sense with reference to the type of calculation of A. Because the first question here would be: is this a determination || determination, or a calculated sentence? Because if 25 $\times$ 25 = 625 is a determination (basic rule), then the multiplication sign now means something different to what it actually means, for example. (I.e.\ we are dealing with a different type of calculation.) And if A is a determination, then this defines the addition differently than if it is a calculated rate. Because the determination is then an explanation of the addition sign \& the calculation rules that allow A to be calculated are a different explanation of the same sign. I must not forget here that $\alpha,\beta,\gamma$ is not the proof of A but only the form of the proof, or of what is proved; $\alpha,\beta,\gamma$ thus defines A.

Therefore, I can only say ``25 $\times$ 25 = 625 is proved" if the method of proof is fixed, independently of the specific proof. Because this method only determines the meaning of ``$\xi\times\eta$", i.e.\ what is proved. In this respect, the form aa : b = c belongs to the method of proof that explains the meaning of c. The question of whether I have calculated correctly is something else. - And so $\alpha,\beta,\gamma$ belong to the method of proof that explains the meaning of theorem A.
\end{quote}
Ms-113,141r [Ms-113: IX, Philosophische Grammatik] (Dated by Wittgenstein from 29.[11.1931] (page 1v) to 23.[5.1932]) (Accessed 30 September 2024)
}
\end{quote}
This was in 1932. But then, in 1941 the correction comes in the following form:
\begin{quote}
Es kommt eben darauf an, was den Sinn des Satzes festlegt. Wovon wir sagen wollen, es lege den Sinn des Satzes fest. Der Gebrauch mu\ss\ ihn festlegen. Aber was rechnen wir zum Gebrauch? -- || Der Gebrauch der Zeichen mu\ss\ ihn festlegen; aber was rechnen wir zum Gebrauch?
 	 			 	
16.6.

     Zwei || Die Beweise beweisen denselben Satz, hei{\ss}t etwa: beide erweisen ihn f\"ur uns als ein brauchbares || geeignetes Instrument zu dem Gleichen. || gleichen Zweck. || als ein geeignetes || passendes Instrument zum gleichen Zweck.
 	 			 	
     Und der Zweck ist eine Anspielung auf Au{\ss}ermathematisches.
 	 			 	
     Ich sagte einmal: `Wenn Du wissen willst, was ein math. Satz sagt, sieh' || schau, was sein Beweis beweist. Nun, ist darin nicht
Wahres \& Falsches? Ist der || Denn ist der Sinn, der Witz, eines math. Satzes wirklich klar, wenn man nur seinen Beweis versteht || sieht? || sieht \& versteht? || , sobald wir nur dem Beweis folgen k\"onnen?\footnote{which can be translated as:
\begin{quote}
It all depends on what determines the meaning of the sentence. What we want to say determines the meaning of the sentence. The use must determine it. But what do we count as usage? - The use of the signs must determine it; but what do we count as use?

16.6
 			 	
     Two proofs prove the same proposition, that is, both prove it to us as a useful instrument for the same purpose. || The same purpose. || as a suitable instrument for the same purpose.
	 			 	
     And the purpose is an allusion to the extra-mathematical.
	 			 	
     I once said: `If you want to know what a math proposition says, see what its proof proves. Well, is'n't it
true \& false? For is the meaning, the point, of a math proposition really clear if one only understands its proof? || sees \& understands? || as soon as we can only follow the proof?
\end{quote}
Ms-124,47 [Ms-124 (WL)] (6 June -- 4 July 1941 and  5 March -- 3 July 1944) (Accessed 30 September 2024)
}
\end{quote}

Wittgenstein's acknowledging the inappropriateness of thinking of proofs and conditions for assertion as determining the meaning of propositions seems to have been overlooked by some `verificationist' theories of meaning which seek a conceptual basis in Wittgenstein's semantic theory and suggest that `meaning is use' is their underlying semantical paradigm. In this context it is important to observe that the emphasis on `verification' in Wittgenstein's moving away from `states of affairs' is rather different from the way those verificationists look at `meaning as verification'. Here is a quote from another manuscript from the 1929-1932 period:
\begin{quote}
Der Sinn eines || des Satzes ist nicht pneumatisch, sondern ist das, was auf die Frage nach der Erkl\"arung des Sinnes zur Antwort kommt. Und -- oder -- der eine Sinn unterscheidet sich vom andern wie die Erkl\"arung des einen von der Erkl\"arung des andern.
	 			 	
 Welche Rolle der Satz im Kalk\"ul spielt, das ist sein Sinn.
 			 	
     Der Sinn steht also nicht hinter ihm (wie der psychische Vorgang der Vorstellungen etc.).
			 	
Welche S\"atze aus ihm folgen \& aus welchen S\"atzen er folgt, das macht seinen Sinn aus. Daher auch die
Frage nach seiner Verifikation eine Frage nach seinem Sinn ist.\footnote{which can be translated by:
\begin{quote}
The sense of a proposition is not pneumatic, but is what comes to the answer to the question of the explanation of the sense. And - or - the one sense differs from the other as the explanation of the one differs from the explanation of the other.
 	 			 	
     What role the proposition plays in the calculation is its sense.
 	 			 	
     The sense is therefore not behind it (like the mental process of ideas etc.).
 	 			 	
Which propositions follow from it \& from which propositions it follows, that constitutes its sense. Hence also the
question about its verification is a question about its meaning.
\end{quote}
Ms-113,42r [Ms-113: IX, Philosophische Grammatik] (Nov 1931 - May 1932) (Accessed 30 September 2024)
}
\end{quote}

Among those proponents of verificationism which take proof conditions (i.e., conditions for assertion) as primitive, is, for example, Michael Dummett (1977, 12):
\begin{quote}
the meaning of each [logical] constant is to be given by specifying, for any sentence in which that constant is the main operator, what is to count as a proof of that sentence, it being assumed that we already know what is to count as a proof of any of the constituents.
\end{quote}
Underlying the so-called Dummett--Prawitz theory of meaning via proofs, D.\ Prawitz (1977), puts forward an interpretation of Wittgenstein's `meaning is use':
\begin{quote}
As pointed out by Dummett, this whole way of arguing with its stress on communication and the role of the language of mathematics is inspired by ideas of Wittgenstein and is very different from Brouwer's rather solipsistic view of mathematics as a languageless activity. Nevertheless, as it seems, it constitutes the best possible argument for some of Brouwer's conclusions. (...)

I have furthermore argued that the rejection of the Platonistic theory of meaning depends, in order to be conclusive, on the development of an adequate theory of meaning along the lines suggested in the above discussion of the principles concerning meaning and use. Even if such a Wittgensteinian theory did not lead to the rejection of classical logic, it would be of great interest in itself.
\end{quote}
Similarly, Dummett (1991, 251f) goes along those lines:
\begin{quote}
Gerhard Gentzen, who, by inventing both natural deduction and the sequent calculus, first taught us how logic should be formalised, gave a hint how to do this, remarking without elaboration that `an introduction rule gives, so to say, a definition of the constant in question', by which he meant that it fixes its meaning, and
that the elimination rule is, in the final analysis, no more than a consequence of this definition. (...) Plainly,
the elimination rules are not consequences of the introduction rules in the straightforward sense of being derivable from them; Gentzen must therefore have
had in mind some more powerful means of drawing consequences. He was also implicitly claiming that
the introduction rules are, in our terminology, self-justifying.
\end{quote}
Right after this, though, Dummett suggests that P.\ Martin-L\"of in his intuitionistic type theory would be advocating a different standpoint, namely, that the rules of elimination, which establish the immediate consequences of a proposition, define the meaning of the logical connective:
\begin{quote}
Intuitively, Gentzen's suggestion that the introduction rules be viewed as fixing the meanings of the logical constants has no more force than the converse suggestion, that they are fixed by the elimination rules; intuitive plausibility oscillates between these opposing suggestions as we move from one logical constant to another. Per Martin-L\"of has, indeed, constructed an entire meaning-theory for the language of mathematics on the basis of the assumption that it is the elimination rules that determine meaning. The underlying idea is that the content of a statement is what you can do with it if you accept it -- what difference learning that it is true will, or at least may, make to you. This is, of course, the guiding idea of a pragmatist meaning-theory. When applied to the logical constants, the immediate consequences of any logically complex statement are taken as drawn by means of an application of one of the relevant elimination rules.
\end{quote}
Taking a similar standpoint, Martin-L\"of himself clearly adheres to the view that the introduction rules, i.e.\ those which establish the assertability conditions, should be seen as the rules defining the meaning of a logical connective:
\begin{quote}
The intuitionists explain the notion of proposition, not by saying that a proposition is the expression of its truth conditions, but rather by saying, in Heyting's words, that a proposition expresses an expectation or an intention, and you may ask, An expectation or an intention of what? The answer is that it is an expectation or an intention of a proof of that proposition. And Kolmogorov phrased essentially the same explanation by saying that a proposition expresses a problem or task (Ger.\ \emph{Aufgabe}). Soon afterwards, there appeared yet another explanation, namely, the one given by Gentzen, who suggested that the introduction rules for the logical constants ought to be considered as so to say the definitions of the constants in question, that is, as what gives the constants in question their meaning. What I would like to make clear is that these four seemingly different explanations actually all amount to the same, that is, they are not only compatible with each other but they are just different ways of phrasing one and the same explanation. (Martin-L\"of (1987, 410.))
\end{quote}
This passage directly recalls the view expressed in his of \emph{Intuitionistic Type Theory} (1984, 24), that the introduction rules define the meaning:
\begin{quote}
The introduction rules say what are the canonical elements (and equal canonical elements) of the set, thus giving its meaning.
\end{quote}

The notion of verification put forward by those proponents of a language-based interpretation of intuitionism is rather different from that of Wittgenstein who in several occasions calls upon the importance of taking seriously the connection between the meaning of a logical symbol and its use/application (Anwendung), as in:
\begin{quote}
In der Logik geschieht immer wieder, was in dem Streit \"uber das Wesen der Definition geschehen ist. Wenn man sagt, die Definition habe es nur mit Zeichen zu tun und ersetze blo\ss\ ein kompliziertes Zeichen durch ein einfacheres -- ein Zeichen durch ein anderes, so wehren sich die Menschen dagegen und sagen, die Definition leiste nicht nur das, oder es gebe eben verschiedene Arten von Definitionen -- der Definition und die interessante und wichtige sei nicht die (reine) ,,Verbaldefinition".

     Sie glauben n\"amlich, man nehme der Definition ihre Bedeutung, Wichtigkeit, wenn man sie als blo{\ss}e Ersetzungsregel, die von Zeichen handelt, hinstellt. W\"ahrend die Bedeutung der Definition in ihrer Anwendung liegt, quasi in ihrer Lebenswichtigkeit. Und eben das geht (heute) in dem Streit zwischen Formalismus, Intuitionismus, etc. vor sich. Es ist den Leuten? unm\"oglich, die Wichtigkeit einer Sache || Handlung || Tatsache, ihre Konsequenzen, ihre Anwendung, von ihr selbst zu unterscheiden; die Beschreibung einer Sache von der Beschreibung ihrer Wichtigkeit.\footnote{which can be translated as:
\begin{quote}
In logic, what happened in the argument about the nature of the definition happens again and again. When it is said that the definition only deals with signs and merely replaces a complicated sign with a simpler one, people object and say that the definition does not only do that, or that there are different kinds of definitions of the definition and that the interesting and important one is not the (pure) ``verbal definition".

     They believe that the definition is deprived of its meaning and importance if it is presented as a mere substitution rule that deals with signs. Whereas the meaning of the definition lies in its application, in its vitality, so to speak. And this is precisely what is going on (today) in the dispute between formalism, intuitionism, etc. It is impossible for people? to distinguish the importance of a thing || action || fact, its consequences, its application, from itself; the description of a thing from the description of its importance.
\end{quote}
Ts-212,XV-108-7 [Ts-212 [cuttings from TSS 208, 210, and 211] (WL)] (1931-1932?)
}
\end{quote}
In a more recent account of his verificationist approach, Martin-L\"of refers to ``Wittgenstein's formula" as the so-called ``verification principle of meaning":
\begin{quote}
It is clear from this what ought to be the general explanation of what a proposition is, namely, that a proposition is defined by stipulating how its proofs, more precisely, canonical or direct proofs, are formed. And, if we take the rules by means of which the canonical proofs are formed to be the introduction rules, I mean, if we call those rules introduction rules as Gentzen did, then his suggestion that the logical constants are defined by their introduction rules is entirely correct, so we may rightly say that a proposition is defined by its introduction rules.

Now what I would like to point out is that this is an explanation which could just as well be identified with the verification principle, provided that it is suitably interpreted. Remember first of all what the verification principle says, namely, that the meaning of a proposition is the method of its verification. The trouble with that principle, considered as a formula, or as a slogan, is that it admits of several different interpretations, so that there arises the question: how is it to be interpreted? Actually, there are at least three natural interpretations of it. On the first of these, the means of verifying a proposition are simply identified with the introduction rules for it, and there is then nothing objectionable about Wittgenstein's formula, provided that we either, as I just did, replace method by means, which is already plural in form, or else make a change in it from the singular to the plural number: the meaning of a proposition is the methods of its verification. Interpreted in this way, it simply coincides with the intuitionistic explanation of what a proposition is, or, if you prefer, the Gentzen version of it in terms of introduction rules. (Martin-L\"of (2013, 6))
\end{quote}
If meant as `meaning is determined by what constitutes a proof', as formalised by Gentzen's `introduction rules' and advocated by Dummett, Prawitz, and Martin-L\"of, this does not take into consideration Wittgenstein's self-corrections.\footnote{In Ms-124 (1941--44) one finds:
\begin{quote}
Ich sagte einmal: `Wenn Du wissen willst, was ein math. Satz sagt, sieh' || schau, was sein Beweis beweist. Nun, ist darin nicht Wahres \& Falsches? Ist der -- Denn ist der Sinn, der Witz, eines math. Satzes wirklich klar, wenn man nur seinen Beweis versteht -- sieht? || sieht \& versteht? -- , sobald wir nur dem Beweis folgen k\"onnen?
\end{quote}
(Ms-124,47 [between 6 June and 4 July 1941 and between 5 March and 3 July 1944] (Accessed 04 October 2024))

\smallskip

\noindent which can be translated as:
\begin{quote}
I once said: ?If you want to know what a math.\ proposition says, see what its proof proves. Well, isn't it
true \& false? For is the meaning, the point, of a math.\ proposition really clear if one only understands its proof? -- sees \& understands? -- as soon as we can only follow the proof? 
\end{quote}
} But here it seems that we need to remember that Wittgenstein's:
\begin{quote}
Die Verifikation ist nicht ein Anzeichen der Wahrheit sondern der Sinn des Satzes.\footnote{which can be translated as:
\begin{quote}
Verification is not a sign of truth but the meaning of the proposition.
\end{quote}
Ms-107,143 [Ms-107: III, Philosophische Betrachtungen] (1929-30) (Accessed 01 October 2024)
The correction is made with respect to the following remark of 1929:
\begin{quote}	
Man k\"onnte fragen: Was sagt (x)2x = x + x? Es sagt da\ss\ alle Gleichungen von der Form 2x = x + x richtig sind. Aber hei{\ss}t das etwas? Kann man sagen: Ja ich sehe da\ss\ alle Gleichungen dieser Form richtig sind, so kann ich jetzt schreiben ,,(x)2x = x + x"?
	 			 	
Ihre Bedeutung mu\ss\ aus ihrem Beweis hervorgehen. Was der Beweis beweist das ist die Bedeutung des Satzes (nicht mehr \& nicht weniger).
\end{quote}
Ms 106,182 [Ms-106: II (\"ONB, Cod.Ser.n.22.019)] (1929) (Accessed 06 October 2024)

which can be translated as:
\begin{quote}
One could ask: What does (x)2x = x + x say? It says that all equations of the form 2x = x + x are correct. But does that mean anything? Can you say: Yes, I see that all equations of this form are correct, so I can now write ``(x)2x = x + x"?
 	 			 	
Their meaning must be clear from their proof. What the proof proves is the meaning of the proposition (no more \& no less).
\end{quote}
}
\end{quote}
is not the same as `proof conditions' or `conditions for assertion', but rather verification as a process which is `external', `from outside', so to speak'. Here is a quote which illustrate this:
\begin{quote}
  Denken wir auch an die Frage ``wie merkst Du, da\ss\ Du Zahnschmerzen hast?", oder gar: ``wie merkst Du, da\ss\ Du f\"urchterliche Zahnschmerzen hast?" (Dagegen: ``wie merkst Du, da\ss\ Du Zahnschmerzen bekommen wirst".)
  
(Hierher geh\"ort die Frage: welchen Sinn hat es, von der Verifikation des Satzes `ich habe Zahnschmerzen' zu reden? Und hier sieht man deutlich, da\ss\ die Frage ``wie wird dieser Satz verifiziert" von einem Gebiet der Grammatik zum andern ihren Sinn \"andert.)

Man k\"onnte nun die Sache so (falsch) auffassen: Die Frage ``wie wei{\ss}t Du, da\ss\ Du Zahnschmerzen hast" wird darum nicht gestellt, weil man dies von den Zahnschmerzen (selbst) aus erster Hand erf\"ahrt, w\"ahrend man, da\ss\ ein Mensch im andern Zimmer ist, aus zweiter Hand, etwa durch ein Ger\"ausch, erf\"ahrt. Das eine wei\ss\ ich durch unmittelbare Beobachtung, das andere erfahre ich indirekt. Also: ``Wie wei{\ss}t Du, da\ss\ Du Zahnschmerzen hast" ? ``Ich wei\ss\ es, weil ich sie habe" -- ``Du entnimmst es daraus, da\ss\ Du sie hast; aber mu{\ss}t Du dazu nicht schon wissen, da\ss\ Du sie hast?". -- -- Der \"Ubergang von den Zahnschmerzen zur Aussage ``ich habe Zahnschmerzen" ist eben ein ganz anderer, als der vom Ger\"ausch zur Aussage ``in diesem Zimmer ist jemand". Das hei{\ss}t, die \"Uberg\"ange geh\"oren ganz andern Sprachspielen an || geh\"oren zu ganz verschiedenen Sprachspielen.

 Ist, da\ss\ ich Zahnschmerzen habe ein Grund zur Annahme, da\ss\ ich Zahnschmerzen habe?\footnote{which can be translated as:
 \begin{quote}
 Let us also think of the question ``How do you realise that you have a toothache?", or even: ``How do you realise that you have a terrible toothache?". (In contrast: ``how do you realise that you will get a toothache").
 
(Here belongs the question: what sense does it make to speak of the verification of the sentence ``I have a toothache"? And here you can clearly see that the question ``how is this sentence verified" changes its meaning from one area of grammar to another).

One could (mis)understand the matter in this way: The question ``how do you know that you have a toothache" is not asked because one learns this from the toothache (itself) at first hand, whereas one learns that a person is in the other room at second hand, for example through a noise. The one I know through direct observation, the other I learn indirectly. Thus: ``How do you know that you have a toothache?" - ``I know it because I have it" - ``You infer that you have it; but must you not already know that you have it?". - The transition from the toothache to the statement ``I have a toothache" is quite different from the transition from the sound to the statement ?there is someone in this room?. In other words, the transitions belong to completely different language games.

 Is the fact that I have a toothache a reason to assume that I have a toothache?
 \end{quote}
 Ts-211,604 [Ts-211 [based on MSS 109, 110, 111, 112, 113 and the beginning of Ms-114] (WL)] (Accessed 01 October 2024)
 }
\end{quote}
Important as it seems to be, the conception which is put forward by the proponents of verificationism (which take conditions for proving as meaning-determining) since Heyting, Gentzen et.\ al saying that an explanation of how to `construct' a proof, or when does one have the conditions for assertion, or even what constitutes a ground for the assertion, shows that it makes it a rather `one-sided' and non-interactive account of meaning and communication through language. And at this point we would like to bring in a quote from the beginning of Ts-211:
\begin{quote}
Es ist n\"amlich die Anschauung aufzugeben, da\ss, um vom Unmittelbaren zu reden, wir von dem Zustand in einem Zeitmoment reden m\"u{\ss}ten. Diese Anschauung ist darin ausgedr\"uckt, wenn man sagt: ``alles, was uns gegeben ist, ist das Gesichtsbild und die Daten der \"ubrigen Sinne, sowie die Erinnerung, in dem gegenw\"artigen Augenblick". Das ist Unsinn; denn was meint man mit dem ``gegenw\"artigen Augenblick"? Dieser Vorstellung liegt vielmehr schon ein physikalisches Bild zu Grunde, n\"amlich das, vom Strom der Erlebnisse, den ich nun in einem Punkt -- an einer Stelle quer durchschneide. Es liegt hier eine \"ahnliche Tendenz und ein \"ahnlicher Fehler vor, wie beim Idealismus (oder Solipsismus).
 	 			 	
     Woher aber diese Tendenz, ``zum Unmittelbaren" kommen zu wollen?
     
     Entspringt sie nicht aus dem Bed\"urfnis, die Verifikation des Satzes verstehen zu wollen, die durch unsere Sprache ganz verschleiert ist.
     
Intuitives Denken, das w\"are so, wie eine Schachpartie auf die Form eines dauernden, gleichbleibenden Zustandes gebracht (ebenso undenkbar).
	 	
     Auf die Frage ``wie hast du das gemeint", k\"onnen eben mehrerlei Antworten kommen:
     
``Ich hab's im Ernst (Spa\ss) gemeint."

``Ich wollte damit sagen; da\ss\ ...  (folgt ein Satz)."

``Ich wollte dich nur aufsitzen lassen."

    Wie geht das vor sich, wenn man einen Satz ausspricht und dabei den anderen nur aufsitzen lassen will? Man spricht, l\"achelt, beobachtet den andern || sieht zu, was der Andere macht, f\"uhlt eine Spannung.
    
     Aber nirgends ist die amorphe Meinung. Diese stellt man sich gleichsam vor, wie den Inhalt eines Tiegels, dessen Aufschrift der Satz ist.
 	 			 	
     ``Inhalt des Satzes."\footnote{which can be translated as:
\begin{quote}
For we must abandon the view that, in order to speak of the immediate, we must speak of the state in a moment of time. This view is expressed when one says: `All that is given to us is the visual image and the data of the other senses, as well as the memory, in the present moment'. This is nonsense, because what do we mean by the `present moment'? Rather, this idea is based on a physical image, namely that of the stream of experiences that I now cut across at one point. There is a similar tendency and a similar error here as with idealism (or solipsism).
 	 			 	
     But where does this tendency to want to get `to the immediate' 'come from?
     
     Does it not arise from the need to want to understand the verification of the proposition, which is completely obscured by our language?
     
Intuitive thinking would be like a game of chess reduced to the form of a permanent, constant state (equally unthinkable).

 The question ?how did you mean that? can be answered in several ways:
 
`I meant it seriously (joking).'

`I meant to say that ... (follows a sentence).'

`I was just trying to get you mounted.'

    What happens when you utter a sentence and just want to let the other person sit up? You speak, smile, watch what the other person does, feel a tension.
    
     But nowhere is the amorphous opinion. You imagine it, as it were, like the contents of a jar labelled with the sentence.

     ``Contents of the sentence."
 \end{quote}
Ts-211,3 (Accessed 02 October 2024)
 }
\end{quote}
For this and other reasons, in a previous opportunity we have expressed our view that constructivist semantics need not and should not be `verificationist', in (at least) Heyting, Gentzen, Prawitz, Dummett's sense. The point to be retained from (e.g.) Brouwer's intuitionism is its constructivism. The 'solipsism' is spurious; but, in this sense, Brouwer and Wittgenstein are correct that understanding must be first person, though anyone's achieving such understanding can of course be within social and worldly contexts, including chalk boards.\footnote{Thanks to an anonymous reviewer of a previous publication.} By focussing exclusively on the method of construction of a proof of a proposition, one misses the information on where and how to challenge the proposition. In a dialogue/game-theooretic  based interpretation of reduction rules for proofs (the elimination rules acting on the result of the application of introduction rules) there is a way in which the elimination rules can be seen as playing the role of the attacker/hearer acting on the result of the act of the defender/speaker. This way one has a formal counterpart to the informal notion of verification but this time two-sided, interactive.

Take, for example, a recent publication by Prawitz (2019b) discussing the notion of `ground for an assertion' in mathematics with a view towards finding ``an explanation of the notion of ground that is independent of the notion of proof":
\begin{quote}
a proof in the intuitionistic tradition is either a construction intended by a proposition, or a realization of such a construction, or a
demonstration of the judgement that a certain constructed object is the construction intended by a proposition. Proofs in the first sense, constructions
intended by propositions, are the terms in which the meanings of propositions are determined. They may be viewed as truth-makers, as G\"oran Sundholm has suggested -- the existence of such an object is what makes the
proposition in question true. The condition for something to be the construction intended by a proposition may then be seen as the truth-condition
of the proposition, not radically di?erent from classical truth-conditions,
as Martin-L\"of [15] has emphasized. Indeed, it can be seen as the constructivization of the corresponding classical truth-condition--while a classical
truth-condition is from a realistic point of view a condition that may be
inaccessible to us, concerning a world independent of us, the intuitionistic
truth-condition concerns something that we may construct and that then
comes in our possession. Thus, the first notion of proof is not primarily an
epistemic one, but is first of all concerned with the meaning of propositions. (...)

The explanation of propositions in terms of constructions intended by
them seems to o?er an account of the notion of ground independent of the
notion of proof, because it seems right to say that when one has realized
or found the intended construction, one is in possession of a ground for
asserting the proposition. (...)

(...) it remains that when we make a conscious inference, we
take ourselves to become justified in asserting the conclusion by having got
a ground for the assertion. Therefore, if a valid inference is to be understood
in the way suggested as involving an operation that yields a ground for the
conclusion when applied to grounds for the premisses, we must be able to
recognize in some way that the operation has this property.\footnote{Prawitz (2019b, 25ff).}
\end{quote}
The offer of an explanation of propositions and of making inferences where the focus relies on `tak[ing] ourselves to become justified in asserting the conclusion by having got a ground for the assertion' gives a different perspective from the examples of Wittgenstein where there is a concern with `one-sided' (\emph{einseitige}) standpoint of language:
\begin{quote}
Unbewu{\ss}te Zahnschmerzen.

Was hei{\ss}t der Satz: ,,ich bin mir meiner Zahnschmerzen bewu{\ss}t".

,,Ich bin mir meiner Armut bewu{\ss}t" $\neq$ ,,ich bin arm". Dagegen:

ich bin mir meiner Zahnschmerzen bewu{\ss}t = ich habe Zahnschmerzen. Es sei denn ich f{\"u}hre eine neue Alternative in meiner Ausdrucksweise ein; dann aber mu\ss\ ich erst ihre Anwendung zeigen sonst habe ich ihr noch keinen Sinn gegeben.
 	 			 	
[zu ,,Schmerzen"]

Mu\ss\ sich denn nicht eine Welt beschreiben lassen, worin der solipsistische Fehler uns weniger nahe liegt. Wo die Tatsachen solche sind, da\ss\ wir weniger leicht zu einer einseitigen Grammatik verf\"uhrt werden?
 	 			 	
In meinen Betrachtungen der Mathematik || \"uber die Mathematik spielen winzige Ver\"anderungen der symbolischen Ausdrucksweise eine Rolle. Was so gesagt || dargestellt klar \& durchsichtig ist, kann, ein wenig anders gesetzt, undurchsichtig oder irref\"uhrend sein.\footnote{which can be translated as:
\begin{quote}
Unconscious toothache.

What does the sentence: ``I am aware of my toothache" mean?

``I am aware of my poverty" $\neq$ ``I am poor". In contrast:

I am aware of my toothache = I have a toothache. Unless I introduce a new alternative in my way of expression; but then I must first show its application otherwise I have not yet given it any meaning.

[on ``pain"]

Must it not be possible to describe a world in which the solipsistic error is less obvious to us? Where the facts are such that we are less easily seduced into a one-sided grammar?
 	 			 	
In my reflections on mathematics, tiny changes in symbolic expression play a role. What is thus said || represented clearly \& transparently, can, placed a little differently, become opaque or
be misleading.
\end{quote}
Ms-114,30r (Ms-114: X, Philosophische Grammatik (WL)) (1932--33) (Accessed 9 Apr 2025)
}
\end{quote}
This specific concern with the one-sided standpoint of language reappears in a manuscript from 1933-36:
\begin{quote}
Die Sprache nun interessiert uns || Wir interessieren uns f\"ur die Sprache, als einen Vorgang nach expliziten Regeln. Denn die philosophischen Probleme sind Mi{\ss}verst\"andnisse, die durch Kl\"arung der Regeln, nach denen wir die Worte gebrauchen wollen, zu beseitigen sind.

Wir betrachten die Sprache von einem einseitigen Standpunkt.\footnote{which can be translated as:
\begin{quote}
Language now interests us || We are interested in language as a process according to explicit rules. This is because philosophical problems are misunderstandings that can be eliminated by clarifying the rules according to which we want to use words.

     We look at language from a one-sided point of view.
 \end{quote}
 Ms-140,24r (Ms-140 [so-called Grosses Format] (WL)) (1933?--36) (Accessed 9 Apr 2025)
 }
\end{quote}
Yet another manuscript from1932--33 brings in this particular concern:
\begin{quote}
Wie versteht man eine Geste? Wenn ich bei irgend einer Gelegenheit sage: ``ich verstehe diese Geste", meine ich da da\ss\ ich sie in Worte oder andere Zeichen \"ubersetzen kann? Nicht Gewi\ss\ nicht immer. Ich charakterisiere ein Erlebnis.
 	 			 	
Das Gesicht ist die Seele des K\"orpers. 
 	 			 	
Der Tonfall der \"Uberzeugung \& die \"Uberzeugung aber auch der Tonfall des Glaubens \& der Glaube \& der  Tonfall der Hoffnung \& die Hoffnung.

Man kann den eigenen Charakter sowenig von Au{\ss}en betrachten || erkennen wie die eigene Schrift. Ich habe zu meiner Schrift eine einseitige Stellung die mich verhindert, sie auf gleichem Fu\ss\ mit anderen Schriften zu sehen \& zu vergleichen.

Wir verzichten auf allgemeine Dogmen \"uber unsern Gegenstand, -- die besonderen Beispiele werfen so viel allgemeines Licht auf ihre Umgebung, als ihnen zukommt.
 	 			 	
``Was ist die richtige Art sein Geld auszugeben?"\footnote{which can be translated as:
\begin{quote}
How do you understand a gesture? When I say on some occasion: `I understand this gesture', do I mean that I can translate it into words or other signs? Not certainly not always. I am characterising an experience.
			 	
The face is the soul of the body.
	 	
The tone of conviction \& conviction but also the tone of faith \& faith \& the  tone of hope \& hope.

You can't see your own character from the outside any more than you can recognise your own writing. I have a one-sided attitude towards my writing that prevents me from seeing \& comparing it on an equal footing with other writings.

We dispense with general dogmas about our subject, -- the particular examples throw as much general light on their surroundings as they deserve.

``What is the right way to spend your money?"
\end{quote}
Ms-156a,49r (Ms-156a [source for Ts-219 and Ms-114 (first part)] (WL)) (1932) (Accessed 8 Apr 2025)
}
\end{quote}.

By contrast to the traditional verificationist point of view, in a later phase, Dummett himself admits the possibility of formulating a theory of meaning for logical constants based on assuming that the elimination rules  give meaning to each constant, as noted by P.\ Schroeder-Heister (2018):
\begin{quote}
Most approaches to proof-theoretic semantics consider introduction rules as basic, meaning giving, or self-justifying, whereas the elimination inferences are justified as valid with respect to the given introduction rules. (...)

One might investigate how far one gets by considering elimination rules rather than introduction rules as a basis of proof-theoretic semantics. Some ideas towards a proof-theoretic semantics based on elimination rather than introduction rules have been sketched by Dummett (1991, Ch.\ 13), albeit in a very rudimentary form.
\end{quote}
An earlier statement of the same line of reasoning is found in:
\begin{quote}
Proof-theoretic semantics has several roots, the most specific one being Gentzen's (1934) remarks that the introduction rules in his calculus of natural deduction define the meanings of logical constants while the elimination rules can be obtained as a consequence of this definition. More broadly, it belongs to the tradition according to which the meaning of a term has to be explained by reference to the way it is used in our language. (Schroeder-Heister (2006, 525))
\end{quote}
These remarks have all attempted to 
insist on justifying the so-called `proof-theoretic semantics' by recalling Wittgenstein's `meaning as use' paradigm. In a recent publication, Prawitz recalls Gentzen's dictum when attempting to define what is a valid argument:
\begin{quote}
Gentzen's idea elaborated: To know the \emph{meaning} of a sentence $A$ in $c$-form is to know that \emph{introductions} of $c$ (that is applications of the introduction rule for $c$) are the \emph{canonical} ways of inferring $A$, which is to say: (1) they are valid inferences, and (2) if $A$ can be proved, it can be proved in that way, that is, by a proof whose final step is an application of the introduction rule for $c$ (which we therefore call a \emph{canonical proof}.). (Prawitz (2019a, 221))
\end{quote}
On the other hand, Martin-L\"of investigates the application of Gentzen's principle\footnote{From Gentzen (1935) (published in von Plato (2017,153):
\begin{quote}
The ``introductions" present, so to say, the ``definitions" of the signs in question, and the ``eliminations" are actually just consequences thereof, expressed more or less as follows: In the elimination of a sign, the proposition the outermost sign of which is in question, must ``be used only as what it means on the basis of the introduction of this sign."
(...)
I think one could show, by making precise this idea, that the E-inferences are, through certain conditions, \emph{unique} consequences of the respective I-inferences.
\end{quote}
} for assessing the ways in which one can give meaning to Frege's assertion sign `$\vdash$':
\begin{quote}
(...) we may already formulate what it is natural to call the correctness condition for assertion, namely the condition under which it is right, and here several terms are possible to use: right, correct, proper. I am going to use them in the same sense. The condition under which it is right, or correct, or proper, to make an assertion is that you know how to perform the task which constitutes the content of the assertion. This is what I have called the correctness condition for assertion in my abstract.

For acts in general it is usually illuminating to ask, What is the purpose of the act? In this case, if we accept the correctness condition that I just gave, What is the purpose of making an assertion? Then we have already to bring in that the speech act involves not only the speaker, but also the hearer, the receiver of the speech act. So, what is it that the assertor wants to achieve, what is the purpose of making an assertion? Well, if we stick to this knowledge account of assertion that I am discussing right now, then the purpose is nothing but to convey to the hearer that the speaker knows how to fulfil the content, the task which makes up the content. The speech act of assertion has no other purpose than to transmit from the speaker to the hearer the information that the speaker knows how to fulfil the task which makes up the content of the assertion, and this succeeds because the speaker must adhere to the correctness condition for assertion that I just formulated. (Martin-L\"of (2019, 229f))
\end{quote}
Note that here the `Peircean' connection between meaning and purpose via a dichotomy between a speaker and a hearer is rather different from the verificationist perspective advocated earlier. After recalling Gentzen's suggestion in:
\begin{quote}
I began by saying that this whole
lecture will be roughly about what the meaning is of the assertion sign. We are used to the fact that when we ask for the meaning of some linguistic construction, it should be visible somehow from the rules that govern that construction, in general Wittgensteinian terms. The first example of this is of course Gentzen's suggestion that the logical operations are defined by their introduction rules. (Martin-L\"of (2019, 234))
\end{quote}
Martin-L\"of concludes by acknowledging that in this case the elimination rules, rather than the introduction rules, are `meaning-determining for the assertion sign':
\begin{quote}
What about the assertion sign? If you did not have this new rule ($C$-elim), you would only have the usual rules of inference, which are of the form ($C$-intro). If you were to take the assertion sign to be determined by these rules, the assertion $\vdash C$ could not mean anything than that $C$ has been demonstrated, has been inferred by the usual inference rules. And that is not how Frege introduced the assertion sign, what Frege meant by the assertion sign. I explained that earlier on: it is the acknowledgement of the truth of a content that the assertion sign expresses. So, we simply cannot explain the assertion sign by referring to the rules governing it if you only have the rules ($C$-intro). But now we are in a better situation, because we also have the rules ($C$-elim), and they are precisely the rules that are meaning-determining for the assertion sign. (id. ibid.)
\end{quote}
Observe here that Martin-L\"of comes close to Peirce's `\emph{Utterer} vs \emph{Interpreter}' by talking about what he calls a `deontological' approach to meaning via a dichotomy `\emph{Speaker} vs \emph{Hearer}':
\begin{quote}
(...) now I come to the commitment account of assertion, which has its origin in Peirce's work during a very early stage of the last century, 1902-03, I think. Peirce's view was that an assertion should be understood as a taking on of responsibility, taking responsibility for the content of the assertion. (id., 230)
\end{quote}
By highlighting the primary importance of an interaction between a speaker and a hearer, Martin-L\"of acknowledges the role of explaining consequences for an `external' perspective: what  the hearer can immediately deduce from what the speaker asserts is central to the definition of meaning. And this is revisited in a more recent account of the assertion sign given by Martin-L\"of in a lecture in 2022:
\begin{quote}
From the contentual point of view, we also have to explain what is the purpose of
uttering an assertion. I take it that to explain the meaning of a complete sentence
is the same as explaining the purpose of the act of uttering it, that is, saying
something by means of it. So the question is, What is the purpose of an utterance
of the assertion $\vdash C$? To answer this question we have to introduce, not only the
speaker, who produces the assertion, but also the hearer, who receives the assertion from the speaker, as indicated in the figure

\centerline{speaker $\longrightarrow$ $\vdash C$ $\longrightarrow$ hearer}

This is necessary, because we cannot explain the meaning, which is to say the
purpose, of an assertion speaking about the speaker alone: the meaning has to do
with the interaction between the speaker and the hearer. The definition of the
purpose that I suggest is that the purpose of the speaker's utterance of $\vdash C$ is to
permit the hearer to request the speaker to do $C$, whereupon the speaker becomes
obligated to do $C$, that is, to fulfil the hearer's request. Now we have, not only
the speaker's act of uttering the assertion, but also the hearer's dual speech act of
requesting the speaker to do $C$.
\end{quote}

At this point we have to remember that, although moving away from using `states of affairs' as the basis for an account for the connection between language and reality, Wittgenstein did not abandon the `language as a calculus' perspective, but acknowledges the primary objective as analysing ordinary language, in spite of its inherent vagueness and imprecisions:
\begin{quote}
(...) in general we don't use language according to strict rules -- -- it hasn't been taught us by means of strict rules, either. We, in our discussions on the other hand, constantly compare language with a calculus proceeding according to exact rules.
 	 			 	
     This is a very one-sided way of looking at language. In practice we very rarely use language as such a calculus. For not only do we not think of the rules of usage -- -- of definitions, etc. -- -- while using language, but when we are asked to give such rules in most cases we aren't able to do so. We are unable clearly to circumscribe the concepts we use; not because we don't know their real definition, but because there is no real ``definition" to them. To suppose that there must be would be like supposing that whenever children play with a ball they play a game according to strict rules.\footnote{Ts-309,41 [so-called Blue Book, von Wright's copy] (1933-34). (Accessed 04 October 2024)}
 \end{quote}
\begin{quote}
When we talk of language as a symbolism used in an exact calculus, that which is in our mind can be found in the sciences and in mathematics. Our ordinary use of language conforms to this standard of exactness only in rare cases. Why then do we in philosophizing constantly compare our use of words with one following exact rules? The answer is that the puzzles which we try to remove always spring from just this attitude towards language.\footnote{Id. Ibid.}
\end{quote}
Numerous passages in the \emph{Nachlass} refer to the conception of language as a rule-based system of calculations which underpins the concept of `explaining the meaning' as in:
\begin{quote}
Wir nennen es ``die Bedeutung des Wortes erkl\"aren" wenn wir es in eine andere Sprache \"ubersetzen, aber auch wenn wir statt seiner eine Geste machen, oder wenn wir auf einen Tr\"ager des Namens weisen; etc.. In soviel verschiedenen Weisen wird der Ausdruck ``Erkl\"arung der Bedeutung" gebraucht.
	 			 	
Wenn man sagt: die Bedeutung eines Wortes sei das, was die Erkl\"arung der Bedeutung erkl\"art,-- so denkt man an diese Erkl\"arung also an das Paradigma eines Schrittes in einem Kalk\"ul. Man denkt sich man k\"onnte sie dem zu erkl\"arenden Zeichen beif\"ugen ja sogar das Zeichen durch sie ersetzen.

     Wenn so die Erkl\"arung mit dem Zeichen (oder statt des Zeichens) wiederholt wird so ist klar da\ss\ sie nicht als, ein f\"ur allemal wirkende, Medizin betrachtet wird (sozusagen als Impfung) sondern als Teil unseres fortlaufenden Kalk\"uls. -- unserer fortlaufenden Kalkulation.
 	 			 	
``Was ist die Bedeutung eines Wortes?" -- ``Was ist der Nutzen eines Gegenstandes?"
	 			 	
Man kann von der Erkl\"arung der Bedeutung sagen da\ss\ sie den Gebrauch des Wortes lehrt.\footnote{which can be translated as:
\begin{quote}
We call it `explaining the meaning of the word' when we translate it into another language, but also when we make a gesture instead of it, or when we point to a bearer of the name; etc.... The expression `explaining the meaning' is used in so many different ways.
	 			 	
When we say that the meaning of a word is that which explains the explanation of the meaning, we think of this explanation as the paradigm of a step in a calculation. One thinks one could add it to the sign to be explained, even replace the sign with it.

     If the explanation is thus repeated with the sign (or instead of the sign), it is clear that it is not regarded as a medicine that works once and for all (as an inoculation, so to speak), but as part of our ongoing calculation. -- of our ongoing calculation.
 	 			 	
``What is the meaning of a word?" - ``What is the use of an object?"
	 			 	
One can say of the explanation of the meaning that it teaches the use of the word.
\end{quote}
Ms-156a,27r [Ms-156a [source for Ts-219 and Ms-114 (first part)] (WL)] (Accessed 18 March 2025)
}
\end{quote}
In a number of manuscripts from the period between 1931 and 1936, the view of language as a calculus is restated but in several opportunities there is reference to the misconception of looking for an ideal language, such as in:
\begin{quote}
Ist der Begriff `rot' undefinierbar? ``Undefinierbar", darunter stellt man sich etwas vor wie unanalysierbar; \& zwar so, als w\"are der betreffende || hier ein Gegenstand unanalysierbar (wie ein chemisches Element). Dann w\"are die Logik also doch eine Art sehr allgemeiner Naturwissenschaft. -- Aber die Unm\"oglichkeit der Analyse entspricht einer von uns angenommenen (festgesetzten) Weise || Art \& Weise der Darstellung.
 	 			 	
     Wir k\"onnten sagen || fragen: ``Wie denn, `undefinierbar'! K\"onnten || K\"onnen wir denn versuchen es zu definieren? Und wie?" --
     
 Es ist von der gr\"o{\ss}ten Bedeutung, da\ss\ wir uns zu einem Kalk\"ul der Logik immer ein Beispiel seiner Anwendung denken, auf welches der Kalk\"ul wirklich eine Anwendung findet, \& da\ss\ wir nicht Beispiele, von denen wir || geben \& sagen, sie seien eigentlich nicht die idealen, diese aber h\"atten wir noch nicht. Das ist das Zeichen einer falschen Auffassung. (Russell \& ich haben, in verschiedener Weise an ihr laboriert. Vergleiche was ich in der ``Abhandlung" || ``Log. phil. Abh." \"uber Elementars\"atze
\& Gegenst\"ande sage.) Kann ich den Kalk\"ul \"uberhaupt verwenden, dann ist dies auch die ideale Verwendung, \& die Verwendung um die es sich handelt. Einerseits will man n\"amlich das Beispiel nicht als das eigentliche anerkennen, weil man in ihm eine Mannigfaltigkeit sieht, der der Kalk\"ul nicht Rechnung tr\"agt. Anderseits ist es doch das Urbild des Kalk\"uls \& er davon hergenommen, \& auf eine getr\"aumte Anwendung kann man nicht warten. Man mu\ss\ sich also eingestehen, welches das eigentliche Urbild || Vorbild des Kalk\"uls ist.

     Nicht aber, als habe man damit einen Fehler begangen, den Kalk\"ul von daher genommen zu haben; sondern der Fehler || . Der Fehler liegt darin, dem Kalk\"ul seine wirkliche || eigentliche Anwendung jetzt nicht zuzugestehen \& sie || , sondern sie f\"ur eine nebulose Ferne || einen idealen Fall zu versprechen.\footnote{which can be translated as:
 \begin{quote}
 Is the term `red' undefinable? `Undefinable' is something like unanalysable, as if the object in question were unanalysable (like a chemical element). Logic would then be a kind of very general natural science after all. - But the impossibility of analysing corresponds to a (fixed) way of representation assumed by us.
 	 			 	
     We could say || ask: `How then, ``undefinable" ! Could || we try to define it? And how?'---
     
 It is of the greatest importance that we always think of an example of its application to a calculus of logic, to which the calculus really finds an application, \& that we do not give examples of which we say they are not actually the ideal ones, but which we do not yet have. This is the sign of a false conception. (Russell \& I have laboured at it in different ways. Compare what I wrote in the `Tractatus' || `Log. phil. Abh.' on elementary propositions
\& objects). If I can use the calculus at all, then this is also the ideal use, \& the use in question. On the one hand, one does not want to recognise the example as the real one, because one sees in it a multiplicity that the calculus does not take into account. On the other hand, it is the archetype of the calculus \& it is taken from it, \& one cannot wait for a dreamed application. One must therefore admit to oneself which is the actual archetype of the calculus.

     Not, however, as if one had made a mistake in having taken the calculus from it; but the mistake is to have taken the calculus from it. The mistake lies in not now conceding to the calculus its real actual application, but in promising it for a nebulous distant ideal case.    
 \end{quote}
 Ms-115,55 [Ms-115: XI, Philosophische Bemerkungen (WL)] (Accessed 19 March 2025)
 }
 
 \end{quote}
The self-correction with respect to the false notion of logical analysis that Russell, Ramsey and himself had conceived where one presumably ``waits for a finite logical analysis of facts, as for a chemical analysis of compounds" had already appeared in an earlier manuscript (Ms-113, 1931--32):
\begin{quote}
Wenn also der Logiker sagt, er habe f\"ur eventuell existierende 6-stellige Relationen in der Arithmetik vorgesorgt, so k\"onnen wir fragen: Was wird denn nun zu dem, was Du vorbereitet hast, hinzukommen || hinzutreten, wenn es seine Anwendung findet || finden wird? Ein neuer Kalk\"ul? -- aber den hast Du ja eben nicht vorbereitet. Oder etwas, was den Kalk\"ul nicht tangiert? -- dann interessiert uns das nicht \& der Kalk\"ul, den Du uns gezeigt hast ist uns Anwendung genug.
	 			 	
Die unrichtige Idee ist, da\ss\ die Anwendung eines Kalk\"uls in der Grammatik der wirklichen Sprache, ihm eine Realit\"at zuordnet, eine Wirklichkeit gibt, die er fr\"uher nicht hatte. || Die unrichtige Idee ist: die Anwendung eines Kalk\"uls auf die wirkliche Sprache
verleihe ihm eine Realit\"at die er fr\"uher || vorher nicht hatte.
	 			 	
Aber, wie gew\"ohnlich in unserem Gebiet, liegt hier der Fehler nicht darin, da\ss\ man etwas Falsches glaubt sondern darin da\ss\ man auf eine irref\"uhrende Analogie hinsieht.
	 			 	
Was geschieht denn, wenn die 6-stellige Relation gefunden wird? Wird quasi ein Metall gefunden, das nun die gew\"unschten (vorher beschriebenen) Eigenschaften (das richtige spezifische Gewicht, die Festigkeit, etc.) hat? Nein; ein Wort wird gefunden, das wir tat\"achlich in unsrer Sprache so verwenden, wie wir etwa den Buchstaben R verwendet haben. ,,Ja, aber dieses Wort hat doch Bedeutung \& ,,R" hatte keine! Wir sehen also jetzt da\ss\ dem ,,R" etwas entsprechen kann". Aber die Bedeutung des Wortes besteht ja nicht darin, da\ss\ ihm etwas entspricht. Au{\ss}er etwa, wo es sich um Namen \& benannten Gegenstand handelt, aber da setzt der Tr\"ager des Namens nur den Kalk\"ul fort, also die Sprache. Und es ist nicht so, wie wenn man sagt: ,,diese Geschichte hat sich tats\"achlich zugetragen, sie war nicht blo{\ss}e Fiktion".
 			 	
Das alles h\"angt auch mit dem falschen Begriff der logischen Analyse zusammen den Russell, Ramsey \& ich hatten. So da\ss\ man auf eine endliche logische Analyse der Tatsachen wartet, wie auf eine chemische von Verbindungen. Eine Analyse, durch die man dann etwa eine 7-stellige Relation wirklich findet, wie ein Element, das tats\"achlich das spezifische Gewicht 7 hat.
	 			 	
Die Grammatik ist f\"ur uns ein reiner Kalk\"ul. (Nicht die Anwendung eines auf die Realit\"at.)\footnote{which can be translated as:
\begin{quote}
So when the logician says that he has provided for any existing 6-digit relations in arithmetic, we can ask: What will be added to what you have prepared when it finds its application? A new calculation? - but you have not prepared it. Or something that does not affect the calculus? - then we are not interested in that \& the calculus you have shown us is application enough for us.
	 			 	
The incorrect idea is that the application of a calculus in the grammar of real language assigns it a reality, gives it a reality that it did not have before. || The incorrect idea is: the application of a calculus to real language
gives it a reality that it did not have before.
	 			 	
But, as usual in our field, the mistake here is not that one believes something false but that one looks at a misleading analogy.

What happens when the 6-digit relation is found? Is a metal found that now has the desired (previously described) properties (the correct specific weight, strength, etc.)? No; a word is found that we actually use in our language in the same way as we used the letter R, for example. `Yes, but this word has meaning \& `R' had none! So now we see that the `R' can correspond to something'. But the meaning of the word is not that something corresponds to it. Except, for example, where it is a matter of names \& named objects, but there the bearer of the name only continues the calculation, i.e. the language. And it is not like saying: `this story actually happened, it was not mere fiction'.
 	 			 	
All this is also related to the false notion of logical analysis that Russell, Ramsey \& I had. So that one waits for a finite logical analysis of facts, as for a chemical analysis of compounds. An analysis by which one then really finds, say, a 7-digit relation, like an element that really has the specific weight 7.
	 			 	
Grammar is pure calculation for us. (Not the application of one to reality.)
\end{quote}
Ms-113,61r [Ms-113: IX, Philosophische Grammatik (\"ONB, Cod.Ser.n.22.022)]
}
\end{quote}
A reference to Carnap's idea of constructing elementary propositions is another moment of self-correction:
\begin{quote}
     Die Idee, Elementars\"atze zu konstruieren (wie dies z.B. Carnap versucht hat) beruht auf einer falschen Auffassung der logischen Analyse. Sie betrachtet das Problem dieser Analyse als das, eine Theorie der Elementars\"atze zu finden. Sie lehnt sich an das an was in der Mechanik || z.B. in der Mechanik geschieht wenn eine Anzahl von Grundgesetzen gefunden wird aus denen das ganze System von S\"atzen || der Mechanik folgt || hervorgeht.
 			 	
Meine eigene Auffassung war falsch: teils, weil ich mir \"uber den Sinn der Worte ,,in einem Satz ist ein logisches Produkt versteckt" (\& \"ahnlicher) nicht klar war, zweitens weil auch ich dachte die logische Analyse m\"usse verborgene Dinge an den Tag bringen
(wie es die chemische \& physikalische tut).\footnote{which can be translated as:
\begin{quote}
     The idea of constructing elementary propositions (as Carnap, for example, attempted to do) is based on a false conception of logical analysis. It regards the problem of this analysis as that of finding a theory of elementary propositions. It is modelled on what happens in mechanics, for example, when a number of fundamental laws are found from which the whole system of propositions of mechanics follows.
	 			 	
My own view was wrong: partly because I was not clear about the meaning of the words `a logical product is hidden in a proposition' (\& similar), secondly because I too thought that logical analysis must bring hidden things to light
(as chemical \& physical analysis does).
\end{quote}
Ms-112,133v [Ms-112: VIII, Bemerkungen zur philosophischen Grammatik (\"ONB, Cod.Ser.n.22.021)
] (Accessed 20 March 2025)
}
\end{quote}

\section{Explanations of consequences as movements within language}
Here the intention is to sketch a picture of a formal counterpart to `meaning is use' on the basis of the idea that explanations of consequences via `movements within language' ought to be taken as a central aspect to Wittgenstein's shift from a picture theory on top of the notion of `state of affairs' and `interpretation of symbols' to the (already subtly present in early writings, as we have shown above)`use of symbols' which underpins his `meaning is use' paradigm. As in the \emph{Investigations} ``every interpretation hangs in the air together with what it interprets, and cannot give it any support. Interpretations by themselves do not determine meaning", as well as in a remark from his transitional period (1929-30):``Perhaps one should say that the expression ``interpretation of symbols" is misleading and one should instead say ``the use of symbols"."  Then, again:
\begin{quote}
Eine Interpretation ist immer nur eine im Gegensatz zu einer andern. Sie h\"angt sich an das Zeichen \& reiht es in ein weiteres System ein.\footnote{which can be translated as:
\begin{quote}An interpretation is always just one in contrast to another. It attaches itself to the sign and incorporates it into another system.
\end{quote}
Ms-110,288 (1930-31) [Ms-110: VI, Philosophische Bemerkungen (WL)]
}
\end{quote}
Thinking of interpretation as translation within language, image or action:
\begin{quote}
Wenn ich sage jedes Bild braucht noch eine Interpretation, so hei{\ss}t Interpretation die \"Ubersetzung in ein weiteres Bild oder in die Tat.\footnote{which can be translated as:
\begin{quote}
If I say every image needs another interpretation, then interpretation means translation into another image or into action.
\end{quote}
Ms-153a,14v [Ms-153a: Anmerkungen [source for MSS 110-112] (WL)] (1931) (Accessed 03 October 2024)
}
\end{quote}
The `movement within language' is called upon in a 1930 manuscript:
\begin{quote}
Zu Grunde liegt allen meinen Betrachtungen (das Gef\"uhl) die Einsicht, da\ss\ der Gedanke einen inneren Zusammenhang mit der Welt hat \& keinen \"au{\ss}eren.

Da\ss\ man also das meint, was man sagt. Hei{\ss}t das aber nicht nur, da\ss\ man sich in der Sprache nicht aus der Sprache, oder in den Gedanken nicht aus den Gedanken herausbewegen kann?
	 			 	
Kann der Zusammenhang zwischen dem ,,Meinen da\ss\ p der Fall ist" \& dem Geschehen
von p noch anders ausgedr\"uckt werden als in der Internen Beziehung eben dieser || jener Ausdr\"ucke? (Ich glaube es nicht.)
			 	
In dem Worte ,,etwas || das \& das meinen" liegt das ganze Problem beschlossen || geschlossen.
 	 			 	
Man f\"uhlt das Stellvertretende an dem Gedanken.
	 			 	
,,Ich dachte, Du w\"urdest zu mir kommen". An diesem Satz mu\ss\ sich alles zeigen lassen was man an exakteren S\"atzen zeigen will.
 	 			 	
In der Sprache wird alles ausgetragen.\footnote{which can be translated as:
\begin{quote}
Underlying all my considerations (the feeling) is the realisation that the thought has an inner connection with the world \& not an outer one.
That one therefore means what one says. But does this not only mean that one cannot move in language out of language, or in thought out of thought?
 	 			 	
Can the connection between `thinking that p is the case' \& the event
of p can be expressed in any other way than in the internal relation of just these || those expressions? (I don't think so.)
 	 			 	
In the word ``something || that \& that mean" lies the whole problem decided || closed.
	 			 	
One feels the vicariousness of the thought.
 	 			 	
`I thought you would come to me'. This sentence must show everything that one wants to show in more exact sentences.
	 	
Everything is carried out in language.
\end{quote}
Ms-108,194 [MS 108 IV. Philosophische Bemerkungen (1930)]
}
\end{quote}
It was around this time when the use of the word `explanation' was brought in to add to the idea of replacing interpretations by explanations as movements within language:
\begin{quote}
Ich will also eigentlich sagen: Es gibt nicht Grammatik \& Interpretation der Zeichen. Sondern soweit von einer Interpretation, also von einer Erkl\"arung der Zeichen, die Rede sein kann, soweit mu\ss\ sie die Grammatik selbst besorgen.

     Denn ich brauchte nur zu fragen: Soll die Interpretation durch S\"atze erfolgen? Und in welchem Verh\"altnis sollen diese S\"atze zu der Sprache stehen die sie schaffen?\footnote{which can be translated as:
\begin{quote}
So what I really want to say is: there is no such thing as grammar and interpretation of signs. Rather, insofar as there can be talk of an interpretation, i.e.\ of an explanation of the signs, it must be provided by the grammar itself.

     For I need only ask: Should the interpretation take place through sentences? And how should these sentences relate to the language they create?
\end{quote}
Ms-109,129 (Wittgenstein Nachlass Ms-109: V, Bemerkungen (WL)) (between 11 August 1930 and 3 February 1931) (Accessed 03 October 2024)
}
\end{quote}
More on `explaining the meaning of the word' appears in the manuscripts of subsequent years (1932-33), such as:
\begin{quote}
Wie lernt man die Bedeutung eines Wortes? Da gibt es viele F\"alle.
 	 			 	
Wir nennen es ``die Bedeutung des Wortes erkl\"aren" wenn wir es in eine andere Sprache \"ubersetzen, aber auch wenn wir statt seiner eine Geste machen, oder wenn wir auf einen Tr\"ager des Namens weisen; etc. In soviel verschiedenen Weisen wird der Ausdruck ``Erkl\"arung der Bedeutung" gebraucht.
 	 			 	
Wenn man sagt: die Bedeutung eines Wortes sei das, was die Erkl\"arung der Bedeutung erkl\"art,-- so denkt man an diese Erkl\"arung also an das Paradigma eines Schrittes in einem Kalk\"ul. Man denkt sich man k\"onnte sie dem zu erkl\"arenden Zeichen beif\"ugen ja sogar das Zeichen durch sie ersetzen.

     Wenn so die Erkl\"arung mit dem Zeichen (oder statt des Zeichens) wiederholt wird so ist klar da\ss\ sie nicht als, ein f\"ur allemal wirkende, Medizin betrachtet wird (sozusagen als Impfung) sondern als Teil unseres fortlaufenden Kalk\"uls. || unserer fortlaufenden Kalkulation.\footnote{which can be translated as:
\begin{quote}
How do you learn the meaning of a word? There are many cases.
 	 			 	
We call it `explaining the meaning of the word' when we translate it into another language, but also when we make a gesture instead of it, or when we point to a bearer of the name; etc.... The expression `explaining the meaning' is used in so many different ways.
 	 			 	
When we say that the meaning of a word is that which explains the explanation of the meaning, we think of this explanation as the paradigm of a step in a calculation. One thinks one could add it to the sign to be explained, even replace the sign with it.

     If the explanation is thus repeated with the sign (or instead of the sign), it is clear that it is not regarded as a medicine that works once and for all (as an inoculation, so to speak), but as part of our ongoing calculation. || of our ongoing calculation.
\end{quote}
Ms-156a,27r [source for Ts-219 and Ms-114 (first part)] (WL)] (not earlier than June 1932 and not later than late autumn 1933)
}
\end{quote}
In another manuscript of around the same time (1933), the so-called \emph{Grosses Format}, one finds a series of remarks suggesting the connection between meaning, use, explanation:
\begin{quote}
Ich will erkl\"aren: der Ort eines Worts in der Grammatik ist seine Bedeutung.
 			 	
Ich kann aber auch sagen: Die Bedeutung eines Wortes ist das, was die Erkl\"arung der Bedeutung erkl\"art.
	 			 	
(``Das was 1 cm$^3$ Wasser wiegt, hat man `1 Gramm' genannt." -- ``Ja was wiegt er denn?")
	 			 	
Die Erkl\"arung der Bedeutung erkl\"art den Gebrauch des Wortes.

     Der Gebrauch des Wortes in der Sprache ist seine Bedeutung.
 	 			 	
Die Grammatik beschreibt den Gebrauch der W\"orter in der Sprache.
	 			 	
Sie verh\"alt sich also zur Sprache \"ahnlich wie die Beschreibung eines Spiels, wie die Spielregeln, zum Spiel.\footnote{which can be translated as:
\begin{quote}
I want to explain: the place of a word in grammar is its meaning.
	 			 	
But I can also say: the meaning of a word is what the explanation of the meaning explains.
	 			 	
(``That which weighs 1 cm$^3$ of water has been called `1 gram'." -- ``Yes, what does it weigh?")
	 			 	
The explanation of the meaning explains the use of the word.

     The use of the word in the language is its meaning.
	 			 	
Grammar describes the use of words in language.
	 			 	
It therefore relates to language in a similar way to the description of a game, like the rules of the game relate to the game.
\end{quote}
Ms-140,15r [Ms-140 [so-called Grosses Format] (WL)] (1933) (Accessed 04 October 2024)
}
\end{quote}
And again in Ms-145 [so-called C1] (1933):
\begin{quote}
 Unsern Standpunkt k\"onnte man kurz so darstellen, da\ss\ f\"ur uns die Bedeutung eines Wortes in der Erkl\"arung liegt die man daf\"ur gegeben hat oder: in der Art wie diese Bedeutung, wie der Gebrauch des Wortes gelernt wurde, \& zwar abgesehen davon welche Wirkung dieses Lernen diese Erkl\"arungen sp\"ater gehabt haben mochten. Wir geben blo\ss\ eine Beschreibung der Formen der Schl\"usselb\"arte ganz abgesehen davon -- \& sehen davon ab ob sie aus dem richtigen Material waren um das Schlo\ss\ \"offnen zu k\"onnen \& davon, welche Wirkungen sie in ihren Schl\"ossern gehabt haben mochten.\footnote{which can be translated as:
 \begin{quote}
 Our point of view could be briefly stated in such a way that for us the meaning of a word lies in the explanation given for it or: in the way in which this meaning, how the use of the word was learnt, apart from what effect this learning of these explanations may have had later. We merely give a description of the shapes of the key beards quite apart from this \& see whether they were made of the right material to open the lock \& what effects they may have had in their locks.
 \end{quote}
 Ms-145,62 Ms-145 [so-called C1] (1933) (Accessed 04 October 2024)
 }
\end{quote}

When it comes to the actual movements within language with a view to explaining the meaning of a proposition, the explanation of which (immediate) consequences appears to be of utmost importance:
\begin{quote}
Der Sinn eines -- des Satzes ist nicht pneumatisch, sondern ist das, was auf die Frage nach der Erkl\"arung des Sinnes zur Antwort kommt. Und -- oder --der eine Sinn unterscheidet sich vom andern wie die Erkl\"arung des einen von der Erkl\"arung des andern.
 			 	
Welche Rolle der Satz im Kalk\"ul spielt, das ist sein Sinn.
 			 	
Der Sinn steht also nicht hinter ihm (wie der psychische Vorgang der Vorstellungen etc.).
	 			 	
Welche S\"atze aus ihm folgen \& aus welchen S\"atzen er folgt, das macht seinen Sinn aus. Daher auch die Frage nach seiner Verifikation eine Frage nach seinem Sinn ist.\footnote{which can be translated as:
\begin{quote}
The sense of a proposition is not pneumatic, but is what comes to the answer to the question of the explanation of the sense. And - or - the one sense differs from the other as the explanation of the one differs from the explanation of the other.

     What role the proposition plays in the calculation is its sense.

     The sense is therefore not behind it (like the mental process of ideas etc.).

Which propositions follow from it \& from which propositions it follows, that constitutes its sense. Hence also the
question about its verification is a question about its meaning.
\end{quote}
Ms-113,42r [Ms-113: IX, Philosophische Grammatik] (Nov 1931 - May 1932) (Accessed 03 October 2024)
}
\end{quote}
Important as it is in itself, this passage also makes a connection with other quotes mentioned earlier regarding the relationship between (method of) verification and meaning.

Another manuscript of the same period carries remarks in the same vein, such as, for example:
\begin{quote}
Man k\"onnte sagen: ``das Hersagen der Regel ist ein Kriterium des Verst\"andnisses, wenn er die Regel mit Verst\"andnis ausspricht \& nicht rein mechanisch." Aber hier kann wieder die sinnvolle Betonung beim Aussprechen als Verst\"andnis gelten; \& warum dann nicht einfach das Aussprechen selbst?

Verstehen = be-greifen = einen bestimmten Eindruck von dem Gegenstand erhalten, ihn auf sich wirken lassen. Einen Satz auf sich wirken lassen; Konsequenzen von ihm || des Satzes betrachten, sich vorstellen; etc..
	 			 	
``Verstehen" nennen wir ein psychisches Ph\"anomen das speziell mit der Erscheinung || den Erscheinungen des Lernens \& Gebrauchs unserer, der menschlichen, Wortsprache verbunden ist.\footnote{which can be translated as:
\begin{quote}
One could say: ?the pronouncing of the rule is a criterion of understanding if he pronounces the rule with understanding \& not purely mechanically.? But here again the meaningful emphasis when pronouncing can count as understanding; \& then why not simply the pronouncing itself?
	 			 	
Understanding = grasping = getting a certain impression of the object, letting it affect you. Let a sentence affect you; consider the consequences of the sentence, imagine it; etc..
	 			 	
We call `understanding' a psychological phenomenon that is specifically connected with the phenomena of learning \& using our, the human, word language.
\end{quote}
Ms-114,123r (1931-32?) [Ms-114: X, Philosophische Grammatik (WL)] (Accessed 03 October 2024)
}
\end{quote}
\begin{quote}
Der Satz sei: ``N ging heute nachmittag ins Senathaus". Der Satz ist f\"ur mich kein blo{\ss}er Laut || Schall || Klang , er ruft in mir eine Vorstellung hervor etwa eines Mannes in der N\"ahe des Senathauses. Aber der Satz \& diese Vorstellung ist nicht blo\ss\ ein Klang \& eine schwache Vorstellung; sondern der Satz hat es sozusagen in sich diese Vorstellung hervorzurufen, aber auch andere Konsequenzen, \& das ist sein Sinn. Die Vorstellung scheint nur ein schwaches Abbild dieses Sinnes, oder, sagen wir, nur eine Ansicht dieses Sinnes. --

Aber was meine ich denn damit; sehe ich den Satz eben nicht als Glied in einem System von Konsequenzen?\footnote{which can be translated as:
\begin{quote}
The sentence is: ``N went to the Senate House this afternoon". The sentence is not a mere sound for me, it evokes in me an idea of a man near the Senate House. But the sentence \& this idea is not merely a sound \& a faint idea; but the sentence has it in itself, so to speak, to evoke this idea, but also other consequences, \& that is its meaning. The idea seems but a faint image of this sense, or, let us say, but a view of this sense.

But what do I mean by this; do I not see the proposition as a link in a system of consequences?
\end{quote}
Ms-114,140v (1931-32?) [Ms-114: X, Philosophische Grammatik (WL)] (Accessed 03 October 2024)
}
\end{quote}
Almost the same remarks appear in Ms-145, but with an important additional paragraph:
\begin{quote}
Kann ich nicht sagen: Wenn ich mir \"uberlege, was war der Sinn eines Satzes -- worin es bestehe da\ss\ der Satz einen Sinn habe, so finde ich da\ss\ es --| er darin liegt da\ss\ so aus dem Satz nach Festsetzungen Konsequenzen hervorgehen, da\ss\ ich den Satz \"ahnlich wie den Zug in einem Spiel betrachte.\footnote{which can be translated as:
\begin{quote}
Can I not say: If I consider what was the meaning of a proposition -- in which it consists that the proposition has a meaning, I find that it lies in the fact that consequences emerge from the proposition according to determinations, that I regard the proposition similarly to the move in a game.
\end{quote}
Ms-145,81 [Ms-145 [so-called C1] (WL)] (``C1 was used while working on the ?first? revision of the Big Typescript (= TS 213 plus handwritten modifications and supplementary material) as begun in MS 114ii. The first page of our notebook is dated 14 October 1933) (Accessed 04 October 2024)
}
\end{quote}
Almost immediately below the quotes just cited one finds a special one, given that it suggests once again the link between movements/transformations between sentences and the explanations of meaning/consequences:
\begin{quote}
``Ich komme am 24ten Dezember nach Wien" || in Wien an", das sind doch nicht blo{\ss}e Worte! Gewi\ss\ nicht; wenn ich sie lese, geht au{\ss}er dem Wahrnehmen der Worte noch Verschiedenes in mir vor sich: ich empfinde etwa Freude, stelle mir etwas vor \& dergleichen. -- Aber ich meine doch nicht blo\ss, da\ss\ mit dem Satz verschiedene mehr oder weniger unwesentliche Begleiterscheinungen vor sich gehn sollen; ich meine, der Satz hat doch einen bestimmten Sinn \& den erfasse ich. || nehme ich wahr. Aber was ist denn dieser bestimmte Sinn? Nun, da\ss\ diese bestimmte Person, die ich kenne, dort \& dort hin kommt, etc.. Ja, \& wenn Du den Sinn angibst, bewegst Du Dich in der grammatischen Umgebung des Satzes umher. Du siehst dann die verschiedenen Transformationen \& Konsequenzen des Satzes als pr\"aformiert an; \& das sind sie, sofern sie in einer Grammatik niedergelegt sind. (Du betrachtest eben den Satz, wie einen Zug eines gegebenen Spiels.)\footnote{which can be translated as:
\begin{quote}
``I will arrive in Vienna on 24 December" are not just words! Certainly not; when I read them, apart from perceiving the words, various other things go on in me: I feel joy, for instance, I imagine something \& the like. -- But I do not merely mean that various more or less insignificant accompanying phenomena are supposed to take place with the sentence; I mean that the sentence has a certain meaning \& that I perceive. || I perceive it. But what is this particular meaning? Well, that this particular person, whom I know, comes there \& there, etc.... Yes, \& when you state the sense, you move around in the grammatical environment of the sentence. You then see the various transformations \& consequences of the sentence as preformed; \& they are, provided they are laid down in a grammar. (You look at the sentence like a move in a given game).
\end{quote}
Ms-114,141r (1931-32?) [Ms-114: X, Philosophische Grammatik (WL)] (Accessed 03 October 2024)
}
\end{quote}
A year later, in Ms-145, there comes again the connection between meaning, understanding, and  `consequences', this time regarding it as characteristic of the thought:
\begin{quote}
Wenn ich dieses Bild statt es durch den Satz hervorrufen zu lassen malte \& es jemandem || dem Andern als Mitteilung statt des Satzes zeigte.

     Ich glaube wenn ich jemandem das oben beschriebene Bild als Mitteilung zeige so wird man wieder geneigt sein zu sagen, er m\"usse es verstehen, \& zwar wird man als Vorgang des Verstehens die \"Ubersetzung in die Wortsprache ansehen. Man wird wieder sagen, das Bild habe einen bestimmten Sinn, es teile einen Gedanken mit \& zwar wird man wieder die Konsequenzen des Bildes als f\"ur den Gedanken charakteristisch erachten.\footnote{which can be translated as:
 \begin{quote}
 If I painted this picture instead of having it evoked by the sentence \& showed it to someone else as a message instead of the sentence.
 
     I believe that if I show someone the picture described above as a message, they will again be inclined to say that they must understand it, \& that the process of understanding will be seen as the translation into the language of words. They will again say that the picture has a certain meaning, that it communicates a thought, and they will again regard the consequences of the picture as characteristic of the thought.
  \end{quote}
  Ms-145,83 [Ms-145 [so-called C1] (WL)] (1933) (Accessed 03 October 2024)
 }
 
 Wenn ich sage der Gedanke, hat etwas Geheimnisvolles, so meine ich nicht das Interessante an den psychischen Ph\"anomenen des Denkens (des Vorstellens, der Phantasie) sondern den Gedanken als scheinbar notwendiges Supplement des Satzes wenn er Sinn hat || des Satzzeichens um ihm Sinn zu geben.
 	 			 	
     Der Gedanke als geheimnisvolles Etwas was die Konsequenzen des Satzes enth\"alt.\footnote{which can be translated as:
\begin{quote}
When I say thought has something mysterious about it, I do not mean what is interesting about the psychic phenomena of thought (imagination, fantasy) but thought as a seemingly necessary supplement to the sentence if it has meaning || of the propositional sign to give it meaning.
	 			 	
     The thought as a mysterious something that contains the consequences of the sentence.
\end{quote}
Ms-145,84  [Ms-145 [so-called C1] (WL)] (1933) (Accessed 03 October 2024)
}

Gem\"a\ss\ den Worten ``ich erfasse den Sinn" oder ``ich denke den Gedanken dieses Satzes" nimmst Du einen Vorgang an der zum Unterschied vom blo{\ss}en Satzzeichen diese Konsequenzen beinhaltet.\footnote{which can be translated as:
\begin{quote}
According to the words ``I grasp the meaning" or ``I think the thought of this sentence" you assume a process which, in contrast to the mere propositional sign, contains these consequences.
\end{quote}
Ms-145,86  [Ms-145 [so-called C1] (WL)] (1933) (Accessed 03 October 2024)
}
\end{quote}
The parallel between the purpose of the sentence and the meaning of language would be given by  transformations of the expression according to the rules, that is what one finds in Ms-146 (1933-34?):
\begin{quote}
Es ist nat\"urlich eine gro{\ss}e Wahrheit darin da\ss\ das System des Gebrauchs der W\"orter ihnen ihre Bedeutung, d.h.\ dem Satz seinen Sinn gibt. Oder auch der Zweck des Satzes -- Sinn der Sprache liegt in den Transformationen des Ausdrucks nach den Regeln. Das Denken eine T\"atigkeit) [Wie der Sinn des Spiels in den Transformationen der Spielstellungen.]\footnote{which can be translated as:
\begin{quote}
It is of course a great truth that the system of using words gives them their meaning, i.e.\ the sentence its sense. Or the purpose of the sentence -- the meaning of language lies in the transformations of the expression according to the rules. Thinking is an activity) [Like the meaning of the game in the transformations of the game positions].
\end{quote}
Ms-146 [so-called C2],34r (1933-34?) (Accessed 03 October 2024)
}
\end{quote}

Next, let us bring in two passages from Ms-132 (1946) which are very telling, both with respect to meaning and explanation of consequences, and with what concerns the idea that a method of verification cannot be `one-sided':
\begin{quote}
``Also steht die Schmerz\"au{\ss}erung
wirklich allein da; da sie durch nichts gerechtfertigt ist?" -- Wenn ich das sage schwebt mir unwillk\"urlich ein Bild vor, das, eines -- des Menschen der eine Schmerz\"au{\ss}erung von sich gibt \& dabei nichts empfindet; kein Wunder, da\ss\ mir ungem\"utlich bei dem Satz zumut ist, die Schmerz\"au{\ss}erung stehe allein da.
 			 	
     Wie wenn ein Wortausdruck ein bestimmtes Bild in uns hervorruft, aber dann f\"ur etwas steht || verwendet wird was dem Bild im normalen Sinn entgegengesetzt ist. Wir werden dann immer wieder vom Wortausdruck auf'ss Bild \& dann wieder vom Wortausdruck auf die tats\"achliche Anwendung blicken \& sagen: ``aber es hei{\ss}t doch das! -- Aber es hei{\ss}t
doch das || das Andere!
 	 			 	
\hfill{15.5.}

     ``Also steht die Schmerz\"au{\ss}erung wirklich allein da; ...?" -- Warum soll ich diese Worte, ``sie || die Schmerz\"au{\ss}erung steht allein da", nicht sagen? Welche Konsequenz haben sie denn? Sie haben ja eben keine Konsequenz.
	 			 	
     Dein -- Mein Spiel bleibt -- bewegt sich ganz in der Sprache.
 	 			 	
     Ich sage mir das Wort ``Schmerz" \& stelle mir den -- einen Schmerz vor; \& sage mir: ``da haben wir doch, was das Wort `Schmerz' bezeichnet--". Gewi\ss, das tue ich. Aber was weiter? -- ; was habe ich damit getan? wozu war es n\"utze? (Ich habe die Schenkungsurkunde an mich ausgefertigt; aber was nun weiter damit?)\footnote{which can be translated as:
 \begin{quote}
 ``So the expression of pain really stands alone, since it is not justified by anything?" - When I say this, an image involuntarily comes to mind, that of a person who makes an utterance of pain \& feels nothing; no wonder I feel uncomfortable with the sentence that the utterance of pain stands alone.
 
Like when a word expression evokes a certain image in us, but then stands for something that is the opposite of the image in the normal sense. We will then look again and again from the word expression to the image \& then again from the word expression to the actual application \
\& say: but it means that! - But it means that || the other!
 	 			 	
\hfill{15.5.}

     ``So the expression of pain really stands alone; ...?" - Why should I not say these words, ``it || the expression of pain stands alone"? What consequence do they have? They have no consequence.
 	 			 	
     Your -- my play remains -- moves entirely in language. 

I say the word `pain' to myself \& imagine the pain; \& say to myself: ``Here we have what the word ``pain" denotes". Certainly, I do that. But what else? - What have I done with it? What was it good for? (I have made out the deed of gift to myself; but what now?)
 \end{quote}
 Ms-121,17r (XVII, Philosophische Bemerkungen) (1938-39?) (Accessed 03 October 2024)
 }
\end{quote}

\begin{quote}
``Wo sp\"urst Du den Kummer?" ? In der Seele. -- Was hei{\ss}t das nur? -- Was f\"ur Konsequenzen ziehen wir aus dieser Ortsbestimmung? || Ortsangabe? Eine ist, da\ss\ wir nicht von einem k\"orperlichen Ort des Kummers reden. Aber wir zeigen || deuten doch auf unsern Leib, als w\"are der Kummer in ihm. Ist das, weil wir ein k\"orperliches Unbehagen sp\"uren? Ich wei\ss\ die Ursache nicht. Aber warum soll ich annehmen, sie sei ein leibliches Unbehagen?\footnote{which can be translated as:
\begin{quote}
``Where do you feel the sorrow?" - In the soul. - What does that mean? - What consequences do we draw from this location? || Where do we find it? One is that we are not talking about a physical place of sorrow. But we point to our body as if the sorrow were in it. Is that because we feel a physical discomfort? I don't know the cause. But why should I assume it is a bodily discomfort?
\end{quote}
Ms-132,64 (1946) (Accessed 03 October 2024)
}

 Es gibt ein Benehmen des Kummers \& Anl\"asse des Kummers. Es gibt auch ein Benehmen \& Anl\"asse, der Hoffnung; der Sehnsucht.
 
``Ich hoffe ..." ist ein Benehmen der Hoffnung. Ist es nun eine Beschreibung meines Seelenzustands? Nun, es l\"a{\ss}t doch auf meinen Seelenzustand schlie{\ss}en. Wenn ich Einem
sage ``Ich hoffe noch immer, er wird kommen" -- so kann der Ander daraus Konsequenzen ziehen, die er etwa so beschreiben kann: ``In seinem gegenw\"artigem Seelenzustand wird er ..."\footnote{which can be translated as:
\begin{quote}
 There is a behaviour of sorrow \& occasions of sorrow. There is also a behaviour, \& occasions of hope; of longing.
 
``I hope ..." is a behaviour of hope. Is it a description of my state of mind? Well, it does indicate my state of mind. If I say to someone
I still hope that he will come? -- the other person can draw consequences from this, which he can describe as follows: ``In his present state of mind, he will come ...".
\end{quote}
Ms-132,89 (1946) (Accessed 03 October 2024)
}
\end{quote}
One can find remarks on the role of drawing consequences in the demonstration of knowing the meaning in a very late manuscript, i.e.\ MS-175, in which Wittgenstein has included dates between 23 September 1950 and 21 March 1951 in the text:
\begin{quote}
Habe ich mich nicht geirrt \& hat nicht Moore vollkommen recht? Habe ich nicht den elementaren Fehler gemacht, || , zu verwechseln, was man denkt, mit dem, was man wei{\ss}? Freilich denke ich nicht ``Die Erde hat einige Zeit vor meiner Geburt schon existiert", aber wei\ss\ ich's drum nicht? Zeige ich nicht, da\ss\ ich's wei{\ss}, indem ich immer die Konsequenzen draus ziehe.
	 			 	
     Wei\ss\ ich nicht auch, da\ss\ von diesem Haus keine Stiege 6 Stock tief in die Erde f\"uhrt, obgleich ich noch nie dran gedacht habe?

  Aber zeigt, da\ss\ ich die Konsequenzen draus ziehe, nicht nur, da\ss\ ich diese Hypothese annehme?\footnote{which can be translated as:
  \begin{quote}
  Wasn't I wrong \& isn't Moore absolutely right? Have I not made the elementary mistake of confusing what one thinks with what one knows? Of course I don't think ?The earth existed some time before I was born?, but isn't that why I know it? Do I not show that I know it by always drawing the consequences from it?
	 			 	
     Do I not also know that there is no staircase leading from this house 6 floors down into the earth, although I have never thought of it?
 	 			 	
     But does it show that I draw consequences from it, not just that I accept this hypothesis?
  \end{quote}
  Ms-175,67v,68r,68v [Ms-175 [so-called Notebook no.1] (WL)] (1950-51) (Accessed 03 October 2024)
  }
\end{quote}

\section{Conclusions}
Having Wittgenstein's \emph{Nachlass} at hand  has proven extremely useful in reexamining not just the unpublished writings, but also the body of published work, to the extent that it has allowed us to identify important aspects of the philosopher's common line of thinking, in spite of various changes of directions. One of those aspects is the association of meaning with use, application, purpose, usefulness of symbols in language. The German terms Gebrauch, Anwendung, Verwendung, Zweck appear in several texts since the \emph{Notebooks} up until the very late manuscripts of 1950--51, including the \emph{Investigations}.\footnote{Here it makes sense to refer to what the Editorial Preface of the revised fourth edition of the \emph{Investigations},
by P. M. S. Hacker and Joachim Schulte says:
\begin{quote}
Anscombe was not consistent in her translation of \emph{Gebrauch}, \emph{Verwendung} and \emph{Anwendung}. We have translated \emph{Gebrauch} by `use', \emph{Verwendung} by `use' or `employment', and \emph{Anwendung} by `application'. `Use' also does service for \emph{ben\"utzen}. In general, however, we have not allowed ourselves to be hidebound by the multiple occurrence of the same German word or phrase in different contexts.
\end{quote}
}
In a sequel to the present essay we shall be addressing this aspect of the trajectory of Wittgenstein's writings.
 As for the findings reported here, apart from raising the initial observations concerning the early signs of Bedeutung, Anwendung, Verwendung, Zweck, it is fair to say that by examining those passages quoted here from his unpublished  as well a published works, strongly support and augment the main thrust of our continuing research on the significance of Wittgenstein's suggestion that meaning is specified by explaining the immediate consequences of a term or statement. 
The findings of such series of reports sound rewarding, at least for the fact that by looking at these new implications from the \emph{Nachlass}, we have helped to further connect Wittgenstein's account of `meaning as use', specified by the `calculus' involved in specifying or using terms (or statements), with the pragmatist approaches to meaning as revealed in Peirce's writings on the interaction between the \emph{Interpreter} and the \emph{Utterer} which seem to bear on the ``game"/``dialogue" approaches to meaning (Lorenzen, Hintikka, Fra\"{\i}ss\'e).\footnote{In this connection, we have also drawn attention to the fact that the so-called \emph{Inversion Principle} (as put forward by Gentzen and Prawitz)  has the same status as ``metalevel" explanation of meaning as do these game/dialogical approaches, including Tarski's ``metalanguage" approach. 
Thus, in Natural Deduction, the reduction rules (which formalise exactly the \emph{Inversion Principle}) should be seen as meaning-giving, not merely as a justification of the elimination rule based on the introduction rule(s).
The reduction rules of labelled natural deduction include the case of the ``Identity types" which were crucial for the more mathematically based recent publications connecting proof theory and homotopy theory. (For the sake of self-containedness, we have added an appendix with some details of such an approach via reduction rules of natural deduction.)

All this is part of an effort to draw from the philosophical works of Wittgenstein an idea that meaning in formal languages  should not be too far apart from  to what the philosopher insists on the connection between use, application, utility of formal symbols in proofs in mathematics. The whole enterprise began almost 4 decades ago in 1987--1988 (de Queiroz (1987),(1988ab)), followed by some more in the 1990s. This view of proofs and meaning led to a proposal of reformulating intuitionistic type theory in the following sense: the introduction of symbols for rewriting paths to the syntax of homotopy type theory was proposed by (de Queiroz and Gabbay (1994), de Oliveira and de Queiroz (1999), de Queiroz and de Oliveira (2011,2014), de Queiroz et al. (2016), Ramos et al. (2017), Ramos et al. (2018), de Veras et al. (2023ab)),  the authors bring in a formal entity named `computational path', first appeared in (de Queiroz and Gabbay (1994)), and show that it can be used to formulate the so-called identity type in an explicit manner. As advocated in (Ramos et al.\ (2021ab), de Veras et al.\ (2023ab)), this allows for making useful bridges between theory of computation, algebraic topology, logic, higher categories, and higher algebra (de Veras et al.\ 2023ab, Mart\'{\i}nez-Rivillas and de Queiroz (2022ab,2023ab)).}

\paragraph*{Acknowledgements.} First and foremost, we acknowledge the role of Marcos Silva (UFPE), the main organiser of \emph{VIII Brazilian Society for Analytic Philosophy Conference}, 22-26 July 2024 https://sites.google.com/view/sbfa-sbpha/olinda-2024\_1, 
Academia Santa Gertrudes, Olinda, Pernambuco, Brazil,
 who kindly presented us with an invitation to give a keynote talk. This served as additional stimulus to write this paper, which, in turn led us to new findings in Wittgenstein's \emph{Nachlass}. 
 Last, but not least, thanks and much appreciation to the anonymous reviewer for the scholarly, careful and competent assessment of the paper.
(Gratefully acknowledged is the free access to automated translation applications such as Google Translator and DeepL.com)

\section*{References}
Dummett, D. 1975. The Philosophical Basis of Intuitionistic Logic. In \emph{Logic Colloquium '73}
Edited by H.E. Rose, J.C. Shepherdson. Studies in Logic and the Foundations of Mathematics
Volume 80, Pages 5--40. Elsevier, https://doi.org/10.1016/S0049-237X(08)71941-4\\
\ \\
Dummett, D. 1977. \emph{Elements of Intuitionism} .Oxford University Press. ISBN: 9780198505242 \\
\ \\
Dummett, M. 1991. \textit{The Logical Basis of Metaphysics}, Harvard University Press, Cambridge (Mass.).\\
\ \\
Gabbay, D.M., de Queiroz, R.J.G.B. 1992. ``Extending the Curry-Howard Interpretation to Linear, Relevant and Other Resource Logics" \emph{Journal of Symbolic Logic} 57(4):1319--1365.\\
\ \\
Gentzen, G. 1935. Untersuchungen \"uber das logische Schlie{\ss}en. I. \emph{Math Z} 39, 176--210. https://doi.org/10.1007/BF01201353. (English translation in J.\ von Plato (2017)).\\
\ \\
Henkin, L. 1961. ``Some remarks on infinitely long formulas". In \emph{Infinitistic Methods (Proc. Sympos. Foundations of Math.)}, Warsaw, pages 167--183.
Pergamon, Oxford.\\
\ \\
Hermes, H. 1959. ``Zum Inversionsprinzip der operativen Logik". In Heyting,
A., editor, \emph{Constructivity in Mathematics}, pages 62--68. North-Holland,
Amsterdam.\\
\ \\
Hilmy, S. 1987. \textit{The Later Wittgenstein: The Emergence of a New Philosophical Method}, Blackwell, Oxford.\\
\ \\
Hintikka,  J. 1968. ``Language-games for quantifiers". In N.\ Rescher, editor, \emph{Studies in Logical Theory}, pages 46--72. Blackwell.\\
\ \\
Hintikka,  J. 1979. ``Quantifiers vs. Quantification Theory". In E. Saarinen (ed) \textit{Game-Theoretical Semantics}. Synthese language library, vol 5. Springer, Dordrecht.\\
\ \\
Hintikka, M. and Hintikka, J. 1986. \textit{Investigating Wittgenstein}.  Basil Blackwell, Oxford.\\
\ \\
Hodges, W. 2001. ``A Sceptical Look". \textit{Proceedings of the Aristotelian Society} Supplementary Volume {75}(1):17--32.\\
\ \\
Howard, W. 1980. ``The formulae-as-types notion of construction".
In H. Curry, J.R. Hindley, J. Seldin (eds.), \emph{To H. B. Curry: Essays on Combinatory Logic, Lambda Calculus, and Formalism}. Academic Press (1980). \\
\ \\
Howard, W. 2014. Wadler's Blog, 2014, https://wadler.blogspot.com/2014/08/howard-on-curry-howard.html\\
\ \\
Keiff, L. 2009. ``Dialogical Logic". In \textit{Stanford Encyclopedia of Philosophy}, Stanford. \linebreak https://plato.stanford.edu/entries/logic-dialogical/\\
\ \\
Krabbe, E. 2001. ``Dialogue Logic Restituted". In W.\ Hodges and E.\ C.\ W.\ Krabbe (eds.) \textit{Dialogue Foundations}, Wiley, pp.\ 33--49.\\
\ \\
Lorenzen, P. 1950. ``Konstruktive Begr\"{u}ndung der Mathematik". \emph{Mathematische Zeitschrift} 53(2):162--202.\\
\ \\
Lorenzen, P. 1955. \textit{Einf\"{u}hrung in die operative Logik und Mathematik}. Die Grundlehren der mathematischen Wissenschaften, vol. 78. Springer-Verlag, Berlin-G\"ottingen-Heidelberg. VII + 298 pp\\
\ \\
Lorenzen, P. 1969. \textit{Normative Logic and Ethics\/}, series B.I-Hochschultaschenb\"ucher. Systematische Philosophie, vol. 236$^*$, Bibliographisches Institut, Mannheim/Z\"urich.\\
\ \\
Martin-L\"of, P. 1984. \emph{Intuitionistic Type Theory}. (Notes taken by G.\ Sambin). Bibliopolis, Napoli. ISBN 978-8870881059.\\
\ \\
Martin-L\"of, P. 1987. ``Truth of a Proposition, Evidence of a Judgement, Validity of a Proof". \textit{Synthese} {73}:407--420.\\
\ \\
Martin-L\"of, P.  2013. ``Verificationism Then and Now". Chapter 1 of M.\ van der Schaar (ed.), \emph{Judgement and the Epistemic Foundation of Logic},  Logic, Epistemology, and the Unity of Science 31, DOI 10.1007/978-94-007-5137-8\_1, Springer 2013.\\
\ \\
Martin-L\"of, P.  2019. ``Logic and Ethics". In Piecha, T.; Schroeder-Heister, P. (ed.), \emph{Proof-Theoretic Semantics: Assessment and Future Perspectives}. Proceedings of the Third T\"ubingen Conference on Proof-Theoretic Semantics, 27--30 March 2019,  Univ T\"ubingen, pp.\ 227--235. http://dx.doi.org/10.15496/publikation-35319.\\
\ \\
Martin-L\"of, P.  2022. ``Corrections of Assertion and Validity of
Inference".  Transcript of a lecture given on 26 October 2022 at the Rolf Schock Symposium in Stockholm.\linebreak https://pml.flu.cas.cz/uploads/PML-Stockholm26Oct22.pdf\\
\ \\
Mart\'{\i}nez-Rivillas, D.O. 2022. \emph{Towards a homotopy domain theory}. PhD thesis,
CIn-UFPE (November 2022). Centro de Inform\'atica, Universidade Federal
de Pernambuco, Recife, Brazil. https://repositorio.ufpe.br/handle/123456789/49221\\
\ \\
Mart\'{\i}nez-Rivillas, D.O., de Queiroz, R.J.G.B. 2022a. ``$\infty$-Groupoid Generated by an Arbitrary Topological $\lambda$-Model". \emph{Logic J.\ of the IGPL} 30(3):465--488. https://doi.org/10.1093/jigpal/jzab015 . Also arXiv:1906.05729\\
\ \\
Mart\'{\i}nez-Rivillas, D.O., de Queiroz, R.J.G.B. 2022b. ``Towards a Homotopy Domain Theory". \emph{Archive for Mathematical Logic} 62:559--579. Nov 2022. https://doi.org/10.1007/s00153-022-00856-0\\
\ \\
Mart\'{\i}nez-Rivillas, D.O., de Queiroz, R.J.G.B. 2023a. ``The Theory of an Arbitrary Higher $\lambda$-Model". \emph{Bulletin of the Section of Logic}, 52(1), 39--58. https://doi.org/10.18778/0138-0680.2023.11 . Also arXiv:2111.07092\\
\ \\
Mart\'{\i}nez-Rivillas, D.O., de Queiroz, R.J.G.B. 2023b. ``Solving homotopy domain equations". (Submitted for publication.) arXiv:2104.01195\\
\ \\
Moore, G. E. 1955. ``Wittgenstein's Lectures in 1930--33". \textit{Mind} {54}(253):1--27.\\
\ \\
de Oliveira, A.G., de Queiroz, R.J.G.B. 1999. ``A Normalization Procedure for the Equational Fragment of
Labelled Natural Deduction". \emph{Logic J.\ of the IGPL} 7(2):173--215.\\
\ \\
de Oliveira, A.G., de Queiroz, R.J.G.B. 2005. ``A new basic set of proof transformations". In S.\ Artemov, H.\ Barringer, A.\ Garcez, L.\ Lamb, \& J.\ Woods (Eds.), \emph{We will show them! Essays in Honour of Dov Gabbay} (Vol. 2, pp.\ 499--528). London: College Publications.\\
\ \\
Pears, D. 1987. \textit{The False Prison.  A Study of the Development of Wittgenstein's Philosophy}. Volume I. Clarendon Press, Oxford.\\
\ \\
Pears, D. 1990. ``Wittgenstein's Holism". \emph{Dialectica} 44(2):165--173.\\
\ \\
Peirce, C. S. 1932. \emph{Collected Papers of Charles Sanders Peirce, Volumes I and II: Principles of Philosophy and Elements of Logic}. C. Hartshorne and
P.\ Weiss (ed.). Harvard Univ Press.\\
\ \\
Pietarinen, A-V. 2014. ``Logical and Linguistic Games from Peirce to Grice to Hintikka". \textit{Teorema}
{33}(2) 121--136.\\
\ \\
von Plato, J. 2017. \emph{Saved from the Cellar.
Gerhard Gentzen's Shorthand Notes on Logic and Foundations of Mathematics.} Springer.\\
\ \\
Prawitz, D. 1965. \textit{Natural Deduction: A Proof-Theoretical Study}. Acta Universitatis Stockholmiensis, Stockholm Studies in Philosophy no.\ 3. Almqvist \& Wiksell, Stockholm, G\"oteborg, and Uppsala.\\
\ \\
Prawitz, D. 1977. ``Meaning and proofs: on the conflict between classical and intuitionistic logic". \textit{Theoria} (Sweden) {XLIII} 2--40.\\
\ \\
Prawitz, D.  2019a. ``Validity of Inferences Reconsidered". In Piecha, T.; Schroeder-Heister, P. (ed.), \emph{Proof-Theoretic Semantics: Assessment and Future Perspectives}. Proceedings of the Third T\"ubingen Conference on Proof-Theoretic Semantics, 27--30 March 2019,  Univ T\"ubingen, pp.\ 213--226. http://dx.doi.org/10.15496/publikation-35319.\\
\ \\
Prawitz, D.  2019b. ``The Fundamental Problem of General Proof Theory". In Piecha, T.; Schroeder-Heister, P. (ed.), \emph{General Proof Theory}. Special Issue of \emph{Studia Logica} 107:11--29. https://doi.org/10.1007/s11225-018-9785-9 \\
\\
de Queiroz, R.J.G.B. 1987. ``Note on Frege's notions of definition and the relationship proof theory vs.\ recursion theory (Extended Abstract)". In \emph{Abstracts of the VIIIth International Congress of Logic, Methodology and Philosophy of Science.} Vol.\ 5, Part I, Institute of Philosophy of the Academy of Sciences of the USSR, Moscow, 1987, pp.\ 69--73.\\
\ \\
de Queiroz, R.J.G.B. 1988a. ``A Proof-Theoretic Account of Programming and the Role of Reduction Rules". \emph{Dialectica} 42(4):265--282.\\
\ \\
de Queiroz, R.J.G.B. 1988b. ``The mathematical language and its semantics: to show the consequences of a proposition is to give its meaning". In P.\ Weingartner, G.\ Schurz, E.\ Leinfellner, R.\ Haller, A.\ H\"ubner, (eds.) \emph{Reports of the thirteenth international Wittgenstein symposium} 18, pp.\ 259--266.\\
\ \\
de Queiroz, R.J.G.B. 1989. ``Meaning, function, purpose, usefulness, consequences -- interconnected concepts (abstract)". In \emph{Abstracts of Fourteenth International Wittgenstein Symposium (Centenary Celebration)}, 1989, p.\ 20. Symposium held in Kirchberg/ Wechsel, August 13--20 1989.\\
\ \\
de Queiroz, R.J.G.B. 1991. ``Meaning as grammar plus consequences". \emph{Dialectica} 45(1):83--86.\\
\ \\
de Queiroz, R.J.G.B. 1992. ``Grundgesetze alongside Begriffsschrift (abstract)", in \emph{Abstracts of Fifteenth International Wittgenstein Symposium}, 1992, pp.\ 15--16. Symposium held in Kirchberg/Wechsel, August 16--23 1992.\\
\ \\
de Queiroz, R.J.G.B. 1994. ``Normalisation and Language Games". \emph{Dialectica} 48(2):83--123.\\
\ \\
de Queiroz, R.J.G.B. 2001. ``Meaning, function, purpose, usefulness, consequences -- interconnected concepts". \emph{Logic J of the IGPL} 9(5):693--734.\\
\ \\
de Queiroz, R.J.G.B. 2008. ``On Reduction Rules, Meaning-as-use, and Proof-theoretic Semantics". \emph{Studia Logica} 90:211--247.\\
\ \\
de Queiroz, R. J. G. B. 2023.``From Tractatus to Later Writings and Back -- New Implications from Wittgenstein's Nachlass" SATS, vol. 24, no.\ 2,  pp.\ 167--203. https://doi.org/10.1515/sats-2022-0016\\
\ \\
de Queiroz, R.J.G.B., Gabbay, D.M. 1994. ``Equality in Labelled Deductive Systems and the functional interpretation of propositional equality". In P.\ Dekker, and M.\ Stokhof, (eds.), \emph{Proceedings of the 9th Amsterdam Colloquium 1994}, ILLC/Department of Philosophy, University of Amsterdam, pp.\ 547--546.\\
\ \\
de Queiroz, R.J.G.B., Gabbay, D.M.  1995. ``The functional interpretation of the existential quantifier"'. \emph{Bull of the IGPL} 3(2--3):243--290.\\
\ \\
de Queiroz, R.J.G.B., Gabbay, D.M.  1997. ``The functional interpretation of modal necessity". In M.\ de Rijke, (ed.), \emph{Advances in Intensional Logic}, Applied Logic Series, Kluwer, 1997, pp.\ 61--91.\\
\ \\
de Queiroz, R.J.G.B., Gabbay, D.M. 1999. Labelled Natural Deduction. In: Ohlbach, H.J., Reyle, U. (eds) \emph{Logic, Language and Reasoning}. Trends in Logic, vol 5. Springer, Dordrecht. https://doi.org/10.1007/978-94-011-4574-9\_10\\
\ \\
de Queiroz, R.J.G.B., Maibaum, T.S.E. 1990. ``Proof Theory and Computer Programming". \emph{Zeitschrift f\"ur mathematische Logik und Grundlagen der Mathematik} 36:389--414. \\
\ \\
de Queiroz, R.J.G.B., Maibaum, T.S.E. 1991. ``Abstract Data Types and Type Theory: Theories as Types". \emph{Zeitschrift f\"ur mathematische Logik und Grundlagen der Mathematik} 37:149--166. \\
\ \\
de Queiroz, R.J.G.B., de Oliveira, A.G. 2011. ``The Functional Interpretation of Direct Computations". \emph{Electronic Notes in Theoretical Computer Science} 269:19--40.\\
\ \\
de Queiroz, R.J.G.B., de Oliveira, A.G. 2014.  ``Natural Deduction for Equality: The Missing Entity". In Pereira, L., Haeusler, E., de Paiva, V. (eds) \emph{Advances in Natural Deduction}. Pages 63--91. Trends in Logic, vol 39. Springer, Dordrecht..\\
\ \\
de Queiroz, R.J.G.B., de Oliveira, A.G., Gabbay, D.M. 2011.  \emph{The Functional Interpretation of Logical Deduction}. Vol. 5 of Advances in Logic series. Imperial College Press / World Scientific, Oct 2011.\\
\ \\
de Queiroz, R.J.G.B., de Oliveira, A.G., Ramos, A.F. 2016. ``Propositional Equality, Identity Types, and Computational Paths". \emph{South Amer. J. of Logic} 2(2):245--296.\\
\ \\
Ramos, A.F. 2018. \emph{Explicit computational paths in type theory}. PhD thesis,
CIn-UFPE (August 2018). Centro de Inform\'atica, Universidade Federal
de Pernambuco, Recife, Brazil. \linebreak https://repositorio.ufpe.br/handle/123456789/32902 (Abstract in: Ramos, A. (2019). Explicit Computational Paths in Type Theory. \emph{Bulletin of Symbolic Logic}, 25(2):213-214. \linebreak doi:10.1017/bsl.2019.2)\\
\ \\
Ramos, A.F., de Queiroz, R.J.G.B., de Oliveira, A.G.  2017. ``On the identity type as the type of computational paths". \emph{Logic J.\ of the IGPL} 25(4):562--584.\\
\ \\
Ramos, A.F., de Queiroz, R.J.G.B., de Oliveira, A.G., de Veras, T.M.L. 2018. ``Explicit Computational Paths". \emph{South Amer.\ J.\ of Logic} 4(2):441--484.\\
\ \\
Ramos, A.F., de Queiroz, R.J.G.B., de Oliveira, A.G. 2021a. ``Computational Paths and the Fundamental Groupoid of a Type". In: \emph{Encontro de Teoria da Computa\c{c}\~ao (ETC)}, 6, Evento Online. Porto Alegre: Sociedade Brasileira de Computa\c{c}\~ao, 2021. p.\ 22--25. ISSN 2595-6116. DOI: https://doi.org/10.5753/etc.2021.16371\\
\ \\
Ramos, A.F., de Queiroz, R.J.G.B., de Oliveira, A.G. 2021b. ``Convers\~ao de Termos, Homotopia, e Estrutura de Grup\'oide". In: \emph{Workshop Brasileiro de L\'ogica (WBL)}, 2, Evento Online. Porto Alegre: Sociedade Brasileira de Computa\c{c}\~ao, 2021. p.\ 33--40. ISSN 2763-8731. DOI: https://doi.org/10.5753/wbl.2021.15776\\
\ \\
Schroeder-Heister, P. 2006. ``Validity concepts in proof-theoretic semantics". \emph{Synthese} 148 (3):525--571.\\
\ \\
Schroeder-Heister, P. 2018. ``Proof-Theoretic Semantics" in \textit{Stanford Encyclopedia of Philosophy}, Stanford University. https://plato.stanford.edu/entries/proof-theoretic-semantics/\\
\ \\
V\"an\"a\"anen, J. 2022. {The Strategic Balance of Games in Logic}. arXiv:2212.01658\\
\ \\
de Veras, T.M.L., Ramos, A.F., de Queiroz, R.J.G.B., de Oliveira, A.G.,  2021. ``Calculation of Fundamental Groups via Computational Paths". In: \emph{Encontro de Teoria da Computa\c{c}\~ao (ETC)}, 6. , 2021, Evento Online. Porto Alegre: Sociedade Brasileira de Computa\c{c}\~ao, 2021 . p.\ 17--21. ISSN 2595-6116. DOI: https://doi.org/10.5753/etc.2021.16370.\\
\ \\
de Veras, T.M.L., Ramos, A.F., de Queiroz, R.J.G.B., de Oliveira, A.G.,  2023a. ``A Topological Application of Labelled Natural Deduction". In \emph{South American J.\ of Logic}. https://www.sa-logic.org/aaccess/ruy.pdf . Also arXiv:1906.09105\\
\ \\
de Veras, T.M.L., Ramos, A.F., de Queiroz, R.J.G.B., de Oliveira, A.G.,  2023b. ``Computational Paths - A Weak Groupoid". \emph{Journal of Logic and Computation}, https://doi.org/10.1093/logcom/exad071 . Also arXiv: 2007.07769.\\
\ \\
Wittgenstein, L. 1974. \textit{Letters to Russell, Keynes and Moore},
Ed. with an Introd. by G.\ H.\ von Wright, (assisted by B.\ F.\ McGuinness), Basil Blackwell, Oxford.\\
\ \\
Wittgenstein, L. 1982. \emph{Last Writings on the Philosophy of Psychology}, vol.\ 1, 1982, vol. 2, 1992, G.H.\ von Wright and H.\ Nyman (eds.), trans. C.G.\ Luckhardt and M.A.E.\ Aue (trans.), Oxford: Blackwell.\\
\ \\
Wittgenstein, L. 1974. \textit{Lectures and Conversations on Aesthetics, Psychology and Religious Belief}, 1966, C.\ Barrett (ed.), Oxford: Blackwell.\\
\ \\
Wittgenstein, L. 1961. \textit{Notebooks 1914--1916,} G.\ H.\ von Wright and G.\ E.\ M.\ Anscombe (eds.), Oxford: Blackwell.\\
\ \\
Wittgenstein, L. 1969. \textit{On Certainty,} G. E. M. Anscombe and G. H. von Wright (eds.), G.E.M.\ Anscombe and D. Paul (trans.), Oxford: Blackwell.\\
\ \\
Wittgenstein, L. 1974. \textit{Philosophical Grammar,} R.\ Rhees (ed.), A.\ Kenny (trans.), Oxford: Blackwell.\\
\ \\
Wittgenstein, L. 1953. \textit{Philosophical Investigations,} G.E.M.\ Anscombe and R.\ Rhees (eds.), G.E.M.\ Anscombe (trans.), Oxford: Blackwell. Revised fourth edition (2009), Translated by
G. E. M. Anscombe, P. M. S. Hacker and Joachim Schulte)\\
\ \\
Wittgenstein, L. (1964) \textit{Philosophical Remarks}, R.\ Rhees (ed.), R.\ Hargreaves and R.\ White (trans.), Oxford: Blackwell.\\
\ \\
Wittgenstein, L. 1971. \textit{ProtoTractatus--An Early Version of Tractatus Logico- Philosophicus}, B.\ F.\ McGuinness, T.\ Nyberg, G.\ H.\ von Wright (eds.), D.\ F.\ Pears and B.\ F.\ McGuinness (trans.), Ithaca: Cornell University Press.\\
\ \\
Wittgenstein, L. 1956. \textit{Remarks on the Foundations of Mathematics}, G.H.\ von Wright, R.\ Rhees and G.\ E.\ M.\ Anscombe (eds.), G.\ E.\ M.\ Anscombe (trans.), Oxford: Blackwell, revised edition 1978.\\
\ \\
Wittgenstein, L. 1980. \textit{Remarks on the Philosophy of Psychology,}, vol.\ 1, G. E. M. Anscombe and G. H. von Wright (eds.), G.\ E.\ M.\ Anscombe (trans.), vol. 2, G. H. von Wright and H. Nyman (eds.), C.\ G.\ Luckhardt and M.A.E.\ Aue (trans.), Oxford: Blackwell.\\
\ \\
Wittgenstein, L. 1922. 
\emph{Tractatus Logico-Philosophicus},  Translated by C.K.\ Ogden.
Kegan Paul, Trench, Trubner \& Co., Ltd. New York: Harcourt, Brace \& Company, Inc.\\
\ \\
Wittgenstein, L. 2016. Wittgenstein, Ludwig: Interactive Dynamic Presentation (IDP) of Ludwig Wittgenstein's philosophical \emph{Nachlass} [wittgensteinonline.no]. Edited by the Wittgenstein Archives at the University of Bergen (WAB) under the direction of Alois Pichler. Bergen: Wittgenstein Archives at the University of Bergen 2016-. http://wab.uib.no/ (accessed 04 October 2024)

\appendix
 
\subsection*{Appendix: Reduction rules and dialogue/games semantics}
The renewed examination of Wittgenstein's \emph{Nachlass} presented here supports and augments the significance of Wittgenstein's suggestion that meaning is determined by explicating immediate consequences of a statement or term. These new findings counter the claim that the so-called Gentzen-Dummett-Prawitz-Martin-L\"of `meaning as determined by assertability conditions / introduction rules' are properly designated as `Wittgensteinian views'.
Instead, these passages strongly suggest connections with the ``pragmatist"/``dialogical" approaches to meaning, already in Peirce's writings on the interaction between the \emph{Interpreter} and the \emph{Utterer}, which likewise appear to bear upon the ``game"/``dialogue" approaches to meaning (Lorenzen, Hintikka, Fra\"{\i}ss\'e). The \emph{Inversion Principle} (advocated by Gentzen and Prawitz) also appears to have the same kind of ``metalevel" explication of meaning as those game/dialogical approaches, including Tarski's ``metalanguage" approach. Hence in Natural Deduction, the reduction rules (which formalise the \emph{Inversion Principle}) should be seen as meaning-giving, not a mere justification of the elimination rule(s) based on the introduction rule(s). The reduction rules of labelled natural deduction extends this approach to `Identity types' which correspond to propositional equality. Significantly, the present examination of these materials make an important step towards a formal counterpart to the ``meaning is use" dictum, while highlighting an important common thread from Wittgenstein's very early to very late writings.

There are at least two well-established non-truth-theoretic semantics dealing with the interface of meaning, knowledge, and logic in the context of dialogues, games, or more generally interaction. One of these is also an alternative perspective on proof theory and meaning theory, advocating that Wittgenstein's ``meaning as use" paradigm, as understood in the context of proof theory, by which the so-called reduction rules (showing the workings of elimination rules on the result of introduction rules) should be seen as appropriate to formalise the explanation of the (immediate) consequences one can draw from a proposition, thus to show the function/purpose/usefulness of its main connective in the calculus of language. To recall D.\ Prawitz' \emph{Natural Deduction} (1965) original formulation: `reduction steps'  are defined as rules which operate on proofs, thus they constitute `metalevel' rules rather than `object level' rules. This suggests they can indeed be seen as playing a `semantical' role, in analogy to several other definitions (Tarski's truth conditions, Lorenzen's dialogical rules, Hintikka's game semantics rules, Abelard-Elo\"{\i}se evaluation game rules.

Previously we have pointed out parallels between our proposal and other approaches to the semantics of logical connectives based on explaining the (immediate) consequences of the corresponding proposition, such as Lorenzen's  (1950, 1955, 1969, p.\ 25) dialogical games and Hintikka's (1968, 1979, p.\ 3) semantical games. The general underlying principle is the logical \emph{Inversion Principle}, uncovered by Gentzen, and later by Lorenzen: the elimination procedure is the exact inverse of the introduction, therefore all that can be asked from an assertion is what is indicated  explaining the elimination procedure. One can look at Lorenzen's dialogical games and try to show how the  \emph{Inversion Principle} resides in the game-approach by comparing it to the reduction-based account of meaning. The correspondence works as follows: (i) the assertion corresponds to the introduction rule(s) (and the respective constructors);
(ii) the attacks correspond to the elimination rule(s) (and the respective destructors);
(iii) the set of defenses correspond to the reduction rules (the effect of elimination rules on the result of introduction rules);
(iv) the relation specifying for each attack the corresponding defense(s) are defined by the result of reduction rules.

Taking from the rules of \emph{reduction} between proofs in a system of Natural Deduction enriched with terms alongside formulas, such as in the so-called Curry-Howard interpretation (Howard (1980)) of which Martin-L\"of's type theory can be seen as an instance, and isolating the terms corresponding the derivations, one can draw the following parallels:
Looking at the conclusion of reduction inference rules, one can take the destructor as being the Attack (or `Nature', its counterpart in the terminology of Hintikka's Game-Theoretical Semantics), and the constructor as being the Assertion (or Hintikka's `Myself'). This way the game-theoretic explanations of logical connectives find direct counterpart in the functional interpretation with the semantics of convertibility:

\medskip

\noindent $\land$-$\beta$-{\it reduction\/}
$$\displaystyle{{\displaystyle{{a:A \qquad b:B} \over
{\langle a,b\rangle:A\land B}}\land\mbox{\it -intr}} \over
{{\tt FST}(\langle a,b\rangle):A}}\land\mbox{\it -elim} \qquad \qquad
\twoheadrightarrow_\beta \qquad \qquad a:A$$
$$\displaystyle{{\displaystyle{{a:A \qquad b:B} \over
{\langle a,b\rangle:A\land B}}\land\mbox{\it -intr}} \over
{{\tt SND}(\langle a,b\rangle):B}}\land\mbox{\it -elim} \qquad \qquad
\twoheadrightarrow_\beta \qquad \qquad b:B$$
Associated rewritings:\\
${\tt FST}(\langle a,b\rangle)=_\beta a$\\
${\tt SND}(\langle a,b\rangle)=_\beta b$\\
\ 
\\
\begin{tabbing}
\underline{Assertion}/{\it Introd.\/} $\quad$ \= \underline{Attack}/{\it Elim.\/} $\qquad\qquad\quad \twoheadrightarrow_\beta \qquad\qquad \qquad
$ \= \underline{Defense} \\
\\
Conjunction (`$\land$'): \\
${A}\land{B}$ \> $\quad \quad$ $L$? \> ${A}$ \\
${A}\land{B}$ \> $\quad \quad$ $R$? \> ${B}$ \\
\\
$\langle a,b\rangle:{A}\land{B}$ \> {\tt FST}$(\langle a,b\rangle)$ $\qquad \twoheadrightarrow_\beta$ \> $a:{A}$ \\
$\langle a,b\rangle:{A}\land{B}$ \> {\tt SND}$(\langle a,b\rangle)$ $\qquad \twoheadrightarrow_\beta$ \> $b:{B}$ \\
\end{tabbing}

\noindent $\lor$-$\beta$-{\it reduction\/}
$$\displaystyle{{\displaystyle{{a:A} \over {{\tt inl}(a):A\lor B}}\lor\mbox{\it -intr} \  \displaystyle{{[x:A]} \atop {f(x):C}} \  \displaystyle{{[y:B]} \atop {g(y):C}}} \over
{{\tt CASE}({\tt inl}(a),\upsilon x.f(x),\upsilon y.g(y)):C}}\lor\mbox{\it -elim} \ \ \twoheadrightarrow_\beta \  \displaystyle{{a:A} \atop {f(a/x):C}}$$
$$\displaystyle{{\displaystyle{{b:B} \over {{\tt inr}(b):A\lor B}}\lor\mbox{\small\it -intr} \  \displaystyle{{[x:A]} \atop {f(x):C}} \  \displaystyle{{[y:B]} \atop {g(y):C}}} \over
{{\tt CASE}({\tt inr}(b),\upsilon x.f(x),\upsilon y.g(y)):C}}\lor\mbox{\small\it -elim} \ \ \twoheadrightarrow_\beta \ \  \displaystyle{{b:B} \atop {g(b/y):C}}$$
Associated rewritings:\\
${\tt CASE}({\tt inl}(a),\upsilon x.f(x),\upsilon y.g(y))=_\beta f(a/x)$\\
${\tt CASE}({\tt inr}(b),\upsilon x.f(x),\upsilon y.g(y))=_\beta g(b/y)$\\
\ 
\\
\begin{tabbing}
\underline{Assertion}/{\it Introd.\/} $\quad$ \= \underline{Attack}/{\it Elim.\/} $\qquad\qquad\quad \twoheadrightarrow_\beta \qquad\qquad \qquad
$ \= \underline{Defense} \\
\\

Disjunction (`$\lor$'): \\
${A}\lor{B}$ \> $\quad \quad$ ? \> ${A}$ \\
${A}\lor{B}$ \> $\quad \quad$ ? \> ${B}$ \\
\\
{\tt inl}$(a):{A}\lor{B}$ \> {\tt CASE}({\tt inl}$(a),\upsilon x.f(x),\upsilon y.g(y))$ $\quad \twoheadrightarrow_\beta$ \> $a:A$, $f(a/x):{C}$ \\
{\tt inr}$(b):{A}\lor{B}$ \> {\tt CASE}({\tt inr}$(b),\upsilon x.f(x),\upsilon y.g(y))$ $\quad \twoheadrightarrow_\beta$ \> $b:B$, $g(b/y):{C}$ \\
\\
\end{tabbing}
Whilst in both {\tt FST}$(\langle a,b\rangle)$ and {\tt SND}$(\langle a,b\rangle)$ the destructors allow access to either of the conjuncts, in {\tt CASE}({\tt inl}$(a),\upsilon x.f(x),\upsilon y.g(y))$ and
 {\tt CASE}({\tt inr}$(b),\upsilon x.f(x),\upsilon y.g(y))$ the destructor is not given access to the either disjunct but must ask for whichever disjunct comes from the introduction.

\medskip

\noindent $\rightarrow$-$\beta$-{\it reduction\/}
$$\displaystyle{{\displaystyle{\ \atop {a:A}} \qquad \displaystyle{{\displaystyle{{[x:A]} \atop {b(x):B}}} \over
{\lambda x.b(x):A\rightarrow B}}\rightarrow\mbox{\it -intr}} \over
{{\tt APP}(\lambda x.b(x),a):B}}\rightarrow\mbox{\it -elim} \quad \twoheadrightarrow_\beta
\quad \displaystyle{{a:A} \atop {b(a/x):B}}$$
Associated rewriting:\\
${\tt APP}(\lambda x.b(x),a)=_\beta b(a/x)$\\
\ 
\\
\begin{tabbing}
\underline{Assertion}/{\it Introd.\/} $\quad$ \= \underline{Attack}/{\it Elim.\/} $\qquad\qquad\quad \twoheadrightarrow_\beta \qquad\qquad \qquad
$ \= \underline{Defense} \\
\\
Implication (`$\rightarrow$'): \\
${A}\rightarrow{B}$ \> $A$ $\quad$ ? \> {$B$} \\
\\
$a:A$\\
$\lambda x.b(x):A\to B$ \> {\tt APP}$(\lambda x.b(x),a)$ $\qquad
\twoheadrightarrow_\beta$ \> $b(a/x):{B}$ \\
\\
\end{tabbing}

\noindent $\forall$-$\beta$-{\it reduction\/}
$$\displaystyle{{\displaystyle{\ \atop {a:D}} \qquad \displaystyle{{\displaystyle{{[x:D]} \atop {f(x):P(x)}}} \over
{\Lambda x.f(x):\forall x^D.P(x)}}\forall\mbox{\it -intr}} \over
{{\tt EXTR}(\Lambda x.f(x),a):P(a)}}\forall\mbox{\it -elim} \quad \twoheadrightarrow_\beta
\quad \displaystyle{{a:D} \atop {f(a/x):P(a)}}$$
Associated rewriting:\\
${\tt EXTR}(\Lambda x.f(x),a)=_\beta f(a/x)$\\
\

\noindent $\exists$-$\beta$-{\it reduction\/}
$$\displaystyle{{\displaystyle{{s:D \quad f(s):P(s)} \over
{\varepsilon x.(f(x),s):\exists x^D.P(x)}}\exists\mbox{\it -intr} \quad
\displaystyle{{[t:D,g(t):P(t)]} \atop {d(g,t):C}}} \over
{{\tt INST}(\varepsilon x.(f(x),s),\sigma g.\sigma t.d(g,t)):C}}\exists\mbox{\it -elim} \  \twoheadrightarrow_\beta\ 
\displaystyle{{s:D,f(s):P(s)} \atop {d(f/g,s/t):C}}$$
Associated rewriting:\\
${\tt INST}(\varepsilon x.(f(x),s),\sigma g.\sigma t.d(g,t))=_\beta d(f/g,s/t)$\\
\

\begin{tabbing}
\underline{Assertion}/{\it Introd.\/} $\quad$ \= \underline{Attack}/{\it Elim.\/} $\qquad\qquad\quad \twoheadrightarrow_\beta \qquad\qquad \qquad
$ \= \underline{Defense} \\
\\
Universal Quantifier (`$\forall$'): \\
$\forall x^{D}.{P}(x)$ \> $s:{D}$ ? \> ${P}(s)$ \\
\\
$s:D$\\
$\Lambda x.f(x):\forall x^{D}.{P}(x)$ \\
\> {\tt EXTR}$(\Lambda x.f(x),s)$ $\qquad
\twoheadrightarrow_\beta$ \> $f(s/x):{P}(s)$ \\
\\
Existential Quantifier (`$\exists$'): \\
$\exists x^{D}.{P}(x)$ \> $\quad \quad$ ? \> $s:{D}, {P}(s)$ \\
\\
\> \> $s:D, f(s):P(s)$\\
$\varepsilon x.(f(x),s):\exists x^{D}.{P}(x)$ \\
\> {\tt INST}$(\varepsilon x.(f(x),s),\sigma g.\sigma t.d(g(t),t))$ $\ \ \twoheadrightarrow_\beta$ \> $d(f(s)/g(t),s/t):{C}$ \\
\\
\end{tabbing}
Whilst in {\tt EXTR}$(\Lambda x.f(x),s)$  the destructor has access to the witness (can use a generic element)
$s$, in {\tt INST}$(\varepsilon x.(f(x),s),\sigma g.\sigma t.d(g(t),t))$ the only option is to eliminate over the witness $s$ which is inside the term $\varepsilon x.(f (x), s)$ built with the constructor because it was chosen by the introduction.

\medskip

\noindent $Id$-$\beta$-{\it reduction\/}
$$\displaystyle{{\displaystyle{{u=_r v:A} \over {r(u,v):Id_A(u,v)}}Id\mbox{\it -intr\/} \qquad
\displaystyle{{[u=_t v:A]} \atop {d(t):C}}} \over
{{\tt REWR}(r(u,v),\sigma t.d(t)):C}}Id\mbox{\it -elim\/} \qquad
\twoheadrightarrow_\beta \qquad
\displaystyle{{u=_r v:A} \atop {d(r/t):C}}$$
Associated rewriting:\\
${\tt REWR}(r(u,v),\sigma t.d(t))=_\beta d(r/t)$\\
\
\begin{tabbing}
\underline{Assertion}/{\it Introd.\/} $\quad$ \= \underline{Attack}/{\it Elim.\/} $\qquad\qquad\quad \twoheadrightarrow_\beta \qquad\qquad \qquad
$ \= \underline{Defense} \\
\\
Propositional Equality (`$Id_A(u,v)$'): \\
$Id_A(u,v)$ \> $\quad \quad$ ? \> $u=_r v:A$ \\
\\
$r(u,v):Id_A(u,v)$ \> {\tt REWR}$(r(u,v),\sigma t.d(t))$ $\ \ \twoheadrightarrow_\beta$ \> $d(r/t):{C}$ \\

\end{tabbing}

In {\tt REWR}$(r(u,v),\sigma t.d(t))$ the destructor has no choice regarding the 
`reason' for $u$ being equal to $v$, since $r$ will have been chosen by the time of the assertion, i.e.\ the application of the introduction rule.

\medskip

Drawing these parallels is intended to exhibit a connection with the ``pragmatist"/``dialogical" approaches to meaning as revealed in Peirce's writings on the interaction between the \emph{Interpreter} and the \emph{Utterer} which seem to pertain to the ``game"/``dialogue" approaches to meaning (Lorenzen, Hintikka, Ehrenfeucht-Fra\"{\i}ss\'e.

\end{document}